\documentstyle[11pt]{article}
\def\BP{\vskip 10pt}
\def\bs{\backslash}

\def\F{\noindent}
\def\newpage{\vfill\eject}
\def\q{\quad\quad}
\def\qq{\quad\quad\quad\quad}

\def\SP{\vskip 5pt}

\def\BP{\vskip 10pt}

\def\split{\begin{array}{ll}}
\def\endsplit{\end{array}}
\def\cen{\centerline}
\def\dm{\displaystyle}
\def\mb{\mbox}
\def\r{\,$\mb}

\def\A{\;\,$}

\def\bB{{\bf B}}
\def\bC{{\bf C}}

\def\bN{{\bf N}}
\def\bNm{{\bf N_-}}
\def\Box{$\slash\slash$}
\def\bR{{\bf R}}

\def\CD{{\cal D}}

\def\Con{C_0^\infty}

\def\f{$\,}
\def\gr{\nabla}

\def\m{\vert}
\def\M{\Vert}
\def\pa{\partial}

\def\Ro{{\bf R}}
\def\Rb{{\bf R}^2}
\def\Rc{{\bf R}^3}

\def\RN{{\bf R}^N}

\def\Sch{Schr\"odinger}

\def\T{\leqno}

\def\tx{{\widetilde x}}

\def\Ag{1}
\def\BDG{2}
\def\BM{3}
\def\DPa{4}
\def\DPb{5}
\def\Ei{6}
\def\IS{7}
\def\JS{9}
\def\JSa{8}
\def\RZ{10}
\def\Sa{11}
\def\Sc{12}
\def\Sd{13}
\def\Wea{14}
\def\Web{15}
\def\Wi{16}
\def\Z{17}

\begin{document}

\centerline{\bf  THE REDUCED WAVE EQUATION IN LAYERED MATERIALS}

\vskip 30pt

\centerline{{\sc Willi J\"{A}GER \ and \ Yoshimi SAIT$\bar { \hbox{o} }$}
       \footnote{The second author was partially supported by Deutche 
       Forschungs Gemeinschaft through SFB 359.}}
 
\vskip 60pt
                                                        
        {\bf 1. \ Introduction}               
\BP

        The mathematical theory of wave in layered media is still posing
     interesting mathematical problems even in the linear, stationary case.
        
        In J\"{a}ger-Sait$\bar{\hbox{o}}$\,[\JS] and [\JSa], we studied 
     the spectrum of the reduced wave operator
$$
          H_0 = -\,\mu_0(x)^{-1}\Delta,                             \T (1.1)
$$
     where \f \mu_0(x) \A is a simple function which takes a two positive
     values \f \mu_{01} \A and \f \mu_{02} \A on \f \Omega_1 \A and 
     \f \Omega_2 $ respectively. Here \f \Omega_{\ell} $, \f \ell = 1, 2 $, 
     are open sets of \f \RN \A such that
$$      
\left\{ \split
        \Omega_1 \cap \Omega_2 = \emptyset, \\
        \overline{\Omega_1} \cup \Omega_2 
                = \Omega_1 \cup \overline{\Omega_2} = \RN,
\endsplit \right.                                       \T (1.2)
$$
     \f \overline{\Omega_{\ell}} \A being the closure of \f \Omega_{\ell} $.
     Under a new condition on the separating surface 
     \f S = \pa\Omega_1 = \pa\Omega_2 $, we have established the limiting 
     absorption principle for \f H_0 \A which implies that \f H_0 \A is
     absolute continuous. Our condition is satisfied, for example, for
     the case where \f S \A is a cylinder.

        In this work we are going to extend the results in [\JS] to the
     multimedia case, the case where \f \mu_0(x) \A can take finitely or
     infinitely many values (see \S2). The limiting absorption principle
     will be established and, again, the operator \f H_0 \A is absolute
     continuous. Also we shall consider short-range or long-range 
     perturbation of \f H_0 $, that is, we shall study the operator
$$
           H = -\,\mu(x)^{-1}\Delta,                        \T (1.3)
$$
     where
$$
         \mu(x) = \mu_0(x) + \mu_1(x)                       \T (1.4)
$$ 
     and \f \mu_1(x) \A is short-range or long-range. In this case we 
     shall prove that the point spectrum, if it exists, is discrete, and 
     the limiting absorption principle holds on any interval which does 
     not contain an eigenvalue.

        As for the study of the reduced wave operators with discontinuous
     coefficients, many works have been done for the stratified media 
     in which the coefficients of the operator are the functions of 
     \f x' \in {\bf R}^k \subset \RN $, \f k < N$.
     Some perturbed operators of the above type have been discussed, too.
     Here we refer Wilcox\,[\Wi], Ben-Artzi-Dermanjian-Guillot\,[\BDG], 
     Weder\,[\Wea], [\Web], DeBi\'{e}vre-Pravica\,[\DPa], [\DPb], 
     Boutet de Monvel-Berthier-Manda [\BM], and Zhang\,[\Z]. 
     In [\DPb] S. DeBi\'{e}vre and D. W. Pravica proved
     there is no point spectrum for the stratified propagators without
     any additional conditions other than sufficient smoothness of the 
     coefficients at infinity.  
     
        It seems that there are rather few results for the nonstratified case.
     Eidus\,[\Ei] was the first to consider the reduced wave operators 
     \f H_0 \A with a {\it cone-shape discontinuity}. He imposed 
     the following assumptions on the separating surface \f S $:
     \ {\it there exist 
     positive constants \f c_1 \A and \f c_2 \A such that
$$
            \big| n_N^{(1)}(x) \big| \ge c_1   \qq (x \in S),       \T (1.5)
$$
     and     
$$
        \m x \cdot n^{(1)}(x) \m \le c_2  \qq (x \in S),            \T (1.6)
$$
     where \f n^{(\ell)}(x) $, \f \ell = 1, 2 $, is the unit outward normal
     of \f \Omega_{\ell} \A at \f x $, and \f x \cdot n^{(1)}(x) \A is the 
     inner product of \f x \A and \f n^{(1)}(x) \A in \f \RN $.}\ Note that 
     a cone having its vertex at the origin and the positive \f x_N$-axis as 
     its axis satisfies (1.5) and (1.6). Under the above assumptions, 
     Eidus\,[\Ei] proved the limiting  absorption principle for \f H_0 $, 
     that is, by denoting by \f R_0(z) \A the resolvent of \f H_0 $, the 
     limits
$$
          \lim_{\eta\downarrow 0} R_0(\lambda \pm i\eta) = R_{0\pm}(\lambda)
                     \qq {\rm in\ \ } \bB(L_{2, 1}(\RN), L_{2, -1}(\RN))
                                                                 \T (1.7)       
$$
     exist for \f \lambda > 0$, where the weighted \f L_2 \A space 
     \f L_{2, t}(\RN)$, \f t \in \bR $, is defined by
$$
      L_{2, t}(\RN) = \{ f \ :\ (1 + \m x \m)^{t}f(x) \in L_2(\RN) \}
                                                                  \T (1.8)         
$$
     with its inner product and norm
$$
\left\{ \split
      {\dm   (f, \ g)_t = \int_{\RN} f(x)\overline{g(x)} 
                                 (1 + |x|)^{2t} \, dx,} \\ *[6pt]
      {\dm   \M f \M_t = (f, \ f)_t^{1/2},} \\ *[6pt]               
\endsplit \right.                                               \T (1.9) 
$$
     and \f \bB(X, Y) \A is the Banach space of all bounded linear operators
     from \f X \A into \f Y $. Then, Sait\={o}\,[\Sd] showed that 
     \f L_{2,1}(\RN) \A and \f L_{2,-1}(\RN) \A in (1.7) can be replaced by     
     \f L_{2, \delta}(\RN) \A and \f L_{2, -\delta}(\RN) \A with 
     \f \delta > 1/2 $, respectively. This means that the limiting 
     absorption principle for \f H_0 \A holds on the same weighted \f L_2 \A 
     spaces as are used for the \Sch\ operator (cf. Agmon\,[\Ag], 
     Ikebe-Sait\={o}\,[\IS] and Sait\={o}\,[\Sa]). Recently Roach-Zhang\,[\RZ]
     has shown that \f u = R^{\pm}(\lambda)f $, where \f \lambda > 0 \A and
     \f f \in L_{2, \delta}(\RN) \A with \f \delta > 1/2$, is characterized 
     as a unique solution of the equation
$$
             (-\mu_0(x)^{-1}\Delta - \lambda)u = f             \T (1.10)
$$
     with the radiation condition
$$
   \lim_{R\to\infty} 
       \frac1R \int_{B_R} |\gr u \mp i\sqrt{\lambda\mu(x)}
                         \tx u|^2\, dx = 0  \qq   (\tx = \frac{x}{|x|}),
                                                                \T (1.11)
$$
     \f B_R \A being the ball with radius \f R \A and center at the origin.
     The condition (1.11) is a natural extension of the radiation condition 
     for the \Sch\ operators ([\IS], [\Sa]). [\RZ] also gave another proof 
     of the limiting absorption principle for \f H_0 $.

        In the recent work [\JS] and [\JSa], we studied the reduced wave 
     operators \f H_0 \A with a {\it cylindrical discontinuity} in which
     the separating surface is assumed to satisfy that    
$$
       (\mu_{02} - \mu_{01})(x \cdot n^{(1)})
          = (\mu_{01} - \mu_{02})(x \cdot n^{(2)}) \ge 0     
                             \q (x \in S).                       \T (1.12)
$$
     The condition (1.12) is satisfied if \f \Omega_1 \A is an infinite 
     cylindrical domain which contains the origin and 
     \f \mu_{02} > \mu_{01} $. Then it has been shown 
     again that \f H_0 \A is absolutely continuous. So far it seems that 
     the absence of the point spectrum can not be obtained without imposing 
     some additional conditions such as (1.5)-(1.6) or (1.12).

        In \S2 we define the reduced wave operator \f H_0 \A with multimedia
     and we state our assumption on the separating surface \f S \A
     and the positive function \f \mu_0 \A in which \f \mu_0 \A can take 
     countably infinite values although the condition is a natural 
     extension of the condition (1.12). \S3 is devoted to showing 
     the limiting absorption principle for the unperturbed operator 
     \f H_0 $. Here the arguments are quite parallel to the one 
     in [8] or [9], and hence we shall omit some of the proof.
     In \S4 we shall discuss the point spectrum of the perturbed 
     operator (1.3). It will be shown that the point spectrum of \f H \A
     is discrete. Also some sufficient conditions for the nonexistence of
     the point spectrum of \f H \A will be given. We shall show in \S5 
     that the limiting absorption principle for \f H \A holds on any 
     closed interval which does not contain the point spectrum.
\SP

        {\sc Acknowledgement.} This work was finished when the second author
     was visiting the University of Heidelberg from October 1994 through 
     March 1995. Here he would like to thank Deutsche Forschungs Gemeinschaft for
     its support through SFB 359.  Also the second author is thankful to 
     Professor Willi J\"{a}ger for his kind hospitality during this period.

$$
\ \ \
$$


        {\bf 2. \ The operators \f H_0 \A and \f H $}               
\BP

        In this section we are going to define a reduced wave operator
$$
             H_0 = -\,\mu_0(x)^{-1}\Delta,                   \T (2.1)
$$
     where \f \mu_0(x) \A is a positive, simple function on \f \RN \A 
     which will be specified below, and its perturbed operator
$$
           H = -\,\mu(x)^{-1}\Delta,                        \T (2.2)
$$
     where
$$
         \mu(x) = \mu_0(x) + \mu_1(x)                       \T (2.3)
$$
     such that \f \mu(x) \A is a positive function on \f \RN \A and 
     \f \mu_1(x) \A decays to \f 0 \A at infinity.

        Let us describe the conditions on \f \mu_0(x) $. Let \f \bN \A be
     all positive integers and let \f \bNm \A be all negative integers. 
     Let \f L \A be a subset of integers satisfying one of the following:
$$
\split
     & (I) \ \ \ \ L = \bNm \cup \{ 0 \} \cup \bN, \\
     & (II) \ \ \, L = \{ L_-, \,  L_- + 1, \, \cdots, \, -1, \, 0 \} 
                                                           \cup \bN, \\
     & (III) \  L = \bNm \cup \{ 0, \, 1, \, \cdots, \, L_+ \}, \\
     & (IV) \ \, L = \{ L_-, \, L_- + 1, \, \cdots, \, -1 \} \cup \{ 0 \} 
                             \cup \{ 1, 2, \cdots, L_+ \}, 
\endsplit                                                   \T (2.4)     
$$
     where \f L_- \in \bNm \cup \{ 0 \} \A and \f L_+ \in \{ 0 \} \cup \bN $.
\BP

        {\bf Assumption 2.1.}  Let \f N \A be a positive integer such 
     that \f N \ge 2 $. Let \f L \A be as in (2.4). For each \f \ell \in L $, 
     let \f \Omega_{\ell} \A be an open set in \f \RN $. Let \f \mu_0 \A 
     be a positive function on \f \RN $. The family 
     \f \{ \Omega_{\ell} \}_{\ell\in L} \A and the function \f \mu_0 \A 
     are assumed to satisfy the following (i) $\sim$ (iii):
\SP
     
        (i) \f \{ \Omega_{\ell} \}_{\ell\in L} \A is a disjoint family of
     open sets of \f \RN \A such that
$$
            \RN = \bigcup_{\ell\in L} \overline{\Omega_{\ell}},    \T (2.5)
$$
     where \f \overline{\Omega_{\ell}} \A is the closure of
     \f \Omega_{\ell} $. For any \f R > 0 $, the open ball \f B_R \A with 
     center at the origin and radius \f R \A is covered by a union
     of a finite number of \f \overline{\Omega_{\ell}} $, i.e., for 
     \f R > 0 \A there is a finite subset \f L_R \A of \f L \A such that
$$
             \Omega_{\ell} \cap B_R = \emptyset \qq (\ell \in L - L_R).
                                                                  \T (2.6)
$$
\SP

        (ii) For each \f \ell \in L $, the boundary \f \pa\Omega_{\ell} \A
     of \f \Omega_{\ell} \A is a disjoint union of two continuous surfaces 
     \f S_{\ell}^{(-)} \A and \f S_{\ell}^{(+)} $, i.e.,
$$
\left\{ \split
     \pa\Omega_{\ell} = S_{\ell}^{(-)} \cup S_{\ell}^{(+)}, \\ *[4pt]
     S_{\ell}^{(-)} \cap S_{\ell}^{(+)} = \emptyset
\endsplit \right.                                              \T (2.7)         
$$
     for \f \ell \in L $, where \f S_{\ell}^{(-)} \A and \f S_{\ell}^{(+)} \A 
     are unions of a finite number of smooth surfaces. Here we assume that 
     \f S_{L_-}^{(-)} = \emptyset \A when \f L \A has the smallest number
     \f L_- \A and \f S_{L_+}^{(+)} = \emptyset \A when 
     \f L \A has the largest number \f L_+ $. Further we assume that
$$ 
       S_{\ell}^{(+)} = S_{\ell+1}^{(-)} \qq (\ell \in L),         \T (2.8)
$$                
     where we set \f S_{\ell+1}^{(-)} = \emptyset \A if 
     \f \ell + 1 \not\in L $.
\SP

        (iii) \f \mu_0 \A is a simple function which takes the value
     \f \nu_{\ell} \A on each \f \Omega_{\ell} $, where \f \nu_{\ell} \A
     is a positive number such that
$$
       0 < m_0 \equiv \inf_{\ell\in L} \nu_{\ell} 
                    \le \sup_{\ell\in L} \nu_{\ell} \equiv M_0 < \infty.
                                                                  \T (2.9)
$$
     Let 
$$
         n^{(\ell)}(x) = ( n_1^{(\ell)}(x), \, n_2^{(\ell)}(x), \,
                      \cdots, \, n_N^{(\ell)}(x)) \q  (\ell \in L),   
                                                              \T (2.10)
$$
     be the unit outward normal of \f \Omega_{\ell} \A at a.e.
     \f x \in \pa\Omega_{\ell} $. Then we assume that
$$
       (\nu_{\ell} - \nu_{\ell+1})(n^{(\ell)}(x) \cdot x) \le 0
                       \q (x \in S_{\ell}^{(+)}, \ \ell \in L)
                                                                \T (2.11)
$$
     although \f \ell \ne L_+ \A if \f L \A has the largest number \f L_+ $,
     where \f n^{(\ell)}(x) \cdot x \A is the usual inner product of
     \f n^{(\ell)}(x) \A and \f x \A in \f \RN $.
\BP

        As for the function \f \mu \A we have
\BP

        {\bf Assumption 2.2.} Let \f \mu \A be a measurable function
     on \f \RN \A satisfying the following (i) and (ii):
\SP

        (i) We have
$$
       0 < \widetilde{m_0} \equiv \inf_{x\in\RN} \mu(x) 
             \le \sup_{x\in\RN} \mu(x) \equiv \widetilde{M_0} < \infty.
                                                                  \T (2.12)
$$
\SP

        (ii) Let \f \mu_1 = \mu - \mu_0 $. Then either \f \mu_1 \A is
     short-range, that is,
$$
       |\mu_1(x)| \le c_1(1 + |x|)^{-1-\epsilon} \qq (x \in \RN),   \T (2.13)
$$
     or  \f \mu_1 \A is long-range, that is, \f \mu_1 \A is differentiable 
     such that  
$$
\left\{ \split
          {\dm |\mu_1(x)| \le c_1(1 + |x|)^{-\epsilon} \qq (x \in \RN), } \\  *[5pt]
          {\dm |\gr \mu_1(x)|
                        \le c_1(1 + |x|)^{-1-\epsilon} \qq (x \in \RN), }
\endsplit \right.                                                  \T (2.14)
$$
     with constants \f c_1, \epsilon > 0 $. Throughout this work we assume 
     that \f 0 < \epsilon < 1/2 \A with no loss of generality.
\BP

        Let \f X_0 \A and  \f X \A be Hilbert spaces given by
$$
\left\{ \split
           X_0 = L_2(\RN; \ \mu_0(x)dx), \\ *[4pt]
           X = L_2(\RN; \ \mu(x)dx).
\endsplit \right.                                            \T (2.15)
$$
     The inner product and norm of \f X_0 \A [\,or \f X $] will be denoted
     by \f ( \ , \ )_{X_0} \A and \f \M \ \M_{X_0} \A [or 
     \f ( \ , \ )_{X} \A and \f \M \ \M_{X} $], respectively. Then define
     the operator \f H_0 \A in \f X_0 \A by
$$
\left\{ \split
        D(H_0) = H^2(\RN), \\ *[4pt]
        H_0u = -\,\mu_0(x)^{-1}\Delta u,
\endsplit \right.                                              \T (2.16)
$$
     where \f D(T) \A is the domain of \f T $, \f H^2(\RN) \A is the second
     order Soblev space on \f \RN $. and \f \Delta u \A is defined
     in the sense of distributions. Similarly the operator \f H \A in 
     \f X \A is given by
$$
\left\{ \split
        D(H) = H^2(\RN), \\ *[4pt]
        Hu = -\,\mu(x)^{-1}\Delta u,
\endsplit \right.                                              \T (2.17)
$$
     Then it is easy to see that \f H_0 \A and \f H \A are selfadjoint
     operators in \f X_0 \A and \f X $, respectively.

        Now we are going to give some examples of 
     \f \{ \Omega_{\ell} \}_{\ell\in L} \A and \f \mu_0 \A which satisfy 
     Assumption 2.1. In the following examples we take \f N = 3 \A although
     the $N$-dimensional versions of these examples can be easily obtained.
\BP

        {\sc Example 2.3.} Let \f L = \bNm \cup \{ 0 \} \cup \bN $. Let 
     \f \{ b_{\ell} \}_{\ell\in \bNm \cup \bN} \A be such that
$$
\left\{ \split
            \cdots < b_m < b_{m+1} < \cdots < b_{-1} < 0 < b_1
                      < \cdots < b_{\ell} \cdots , \\ *[4pt]
            b_{\ell} \to \infty \q \ \ \ \ (\ell \to \infty), \\ *[4pt]               
            b_{m} \to -\infty \q (m \to -\infty),
\endsplit \right.                                               \T (2.18) 
$$
     and define \f \{ \Omega_{\ell} \}_{\ell\in L} \A by
$$
\left\{ \split
      \Omega_{\ell} = \{ x = (x_1, x_2, x_3) \in \Rc : 
             b_{\ell} < x_3 <  b_{\ell+1} \} \q (\ell \in \bN), \\ *[4pt]
      \Omega_0 = \{ x = (x_1, x_2, x_3) \in \Rc : 
             b_{-1} < x_3 <  b_{1} \}, \\ *[4pt]
      \Omega_{\ell} = \{ x = (x_1, x_2, x_3) \in \Rc : 
             b_{\ell-1} < x_3 <  b_{\ell} \} \q (\ell \in \bNm), 
\endsplit \right.                                               \T (2.19)
$$
     Then the separating surfaces \f S_{\ell}^{(\pm)} \A are given by
$$
\left\{ \split
      S_{\ell}^{(+)} = \{ x = (x_1, x_2, x_3) \in \Rc : 
              x_3 = b_{\ell+1} \}, \\ *[4pt]
      S_{\ell}^{(-)} = \{ x = (x_1, x_2, x_3) \in \Rc : 
              x_3 =  b_{\ell} \}
\endsplit \right.                                               \T (2.20)
$$     
     for \f \ell \in \bN $,

$$
\left\{ \split
      S_{0}^{(+)} = \{ x = (x_1, x_2, x_3) \in \Rc : 
              x_3 = b_{1} \}, \\ *[4pt]
      S_{0}^{(-)} = \{ x = (x_1, x_2, x_3) \in \Rc : 
              x_3 =  b_{-1} \},
\endsplit \right.                                               \T (2.21)
$$     
     and
$$
\left\{ \split
      S_{\ell}^{(+)} = \{ x = (x_1, x_2, x_3) \in \Rc : 
              x_3 = b_{\ell} \}, \\ *[4pt]
      S_{\ell}^{(-)} = \{ x = (x_1, x_2, x_3) \in \Rc : 
              x_3 =  b_{\ell-1} \}
\endsplit \right.                                               \T (2.22)
$$    
     for \f \ell \in \bNm $. Define \f \mu_0 \A by
$$
          \mu_0(x) = \nu_{\ell} \qq (x \in \Omega_{\ell})       \T (2.23)
$$
     such that 
$$
\left\{ \split
        \nu_0 < \nu_1 < \nu_2 < \cdots < \nu_{\ell} < \cdots, \\ *[4pt]
        \nu_0 < \nu_{-1} < \nu_{-2} < \cdots < \nu_{-\ell} < \cdots
\endsplit \right.                                              \T (2.24)
$$
     with
$$
\left\{ \split
        \nu_0 > 0, \\ *[4pt]
        {\dm M_0 \equiv \sup_{\ell\in\bN\cup\bNm} \nu_{\ell} < \infty. }
\endsplit \right.                                              \T (2.25)
$$
     Since we have
$$
\left\{ \split
        n^{(\ell)}(x)\cdot x \ge 0 \q (x \in S_{\ell}^{(+)}, \ 
                                  \ell \in \{ 0 \} \cup \bN) \\ *[4pt]
        n^{(\ell)}(x)\cdot x \le 0 \q (x \in S_{\ell}^{(+)}, \ 
                                  \ell \in \bNm),
\endsplit \right.                                              \T (2.26)
$$
     we see that the condition (2.11) is satisfied. Although this is a 
     reduced wave operator in stratified media studied by many authors
     (see, e.g., [\Wi], [\Wea], [\DPa]), note that \f \Omega_{\ell} \A 
     can be modified as far as the condition (2.26) holds good.
\BP

        {\sc Example 2.4.}  Let \f L = \{ 0 \} \cup \bN $. Let 
     \f \{ b_{\ell} \}_{\ell\in \bN} \A be such that
$$
\left\{ \split
            0 < b_1 < b_2 < \cdots < b_{\ell} < \cdots, \\ *[4pt]
            b_{\ell} \to \infty \q (\ell \to \infty),                
\endsplit \right.                                               \T (2.27) 
$$
     and define \f \{ \Omega_{\ell} \}_{\ell\in L} \A by
$$
\left\{ \split
      \Omega_{0} = \{ x = (x_1, x_2, x_3) \in \Rc : 
              x_1^2 + x_2^2  <  b_{1}^2 \}, \\ *[4pt]
      \Omega_{\ell} = \{ x = (x_1, x_2, x_3) \in \Rc : 
             b_{\ell}^2 < x_1^2 + x_2^2 < b_{\ell+1}^2 \}, 
                                 \q (\ell \in \bN)                   
\endsplit \right.                                               \T (2.28)
$$
     The separating surfaces \f S_{\ell}^{(\pm)} \A are given by
$$
\left\{ \split
      S_{\ell}^{(+)} = \{ x = (x_1, x_2, x_3) \in \Rc : 
              x_1^2 + x_2^2 = b_{\ell+1}^2 \} \q (\ell \in L), \\ *[4pt]
      S_{\ell}^{(-)} = \{ x = (x_1, x_2, x_3) \in \Rc : 
              x_1^2 + x_2^2 = b_{\ell}^2 \} \q (\ell \in \bN), \\ *[4pt]
      S_{0}^{(-)} = \emptyset,
\endsplit \right.                                               \T (2.29)
$$     
     Define \f \mu_0 \A by
$$
          \mu_0(x) = \nu_{\ell} \qq (x \in \Omega_{\ell})       \T (2.30)
$$
     such that 
$$
        \nu_0 < \nu_1 < \nu_2 < \cdots < \nu_{\ell} < \cdots,  \T (2.31)
$$
     with
$$
\left\{ \split
        \nu_0 > 0, \\ *[4pt]
        {\dm M_0 \equiv \sup_{\ell\in\bN} \nu_{\ell} < \infty. }
\endsplit \right.                                              \T (2.32)
$$     
     Since     
$$
        n^{(\ell)}(x)\cdot x \ge 0 \q (x \in S_{\ell}^{(+)}, \ 
                                  \ell \in L),   \T (2.33)
$$
     from (2.31) it is seen that the condition (2.11) is satisfied. 
     Again \f \Omega_{\ell} \A are allowed to be deformed as far as
     (2.33) holds good.

$$
\ \ \ \
$$                                                          

        {\bf 3. \ The unperturbed operator \f H_0 $}              
\BP

        In this section we are going to discuss the unperturbed operator 
     \f H_0 \A given by (2.16). First we shall show the uniqueness theorem 
     for the equation 
$$
          (-\,\mu_0(x)^{-1}\Delta - \lambda)u = f                   \T (3.1)                                               
$$
     with radiation condition. Then, after showing several a priori estimates 
     of the solution \f u \A of the equation (3.1), the limiting absorption 
     principle for \f H_0 \A will be proved. The arguments in this section
     are quite parallel to the ones in J\"{a}ger-Sait$\bar{\hbox{o}}$\, 
     [\JS], and hence we shall omit the proof or give a sketch of
     proof in most of the theorems given in this section.

        We shall start with some notations.
\BP

        {\sc Notation 3.1.}\ \ Let \f z \in \bC $, \f x = (x_1, x_2, \cdots,
     x_N) $, \f r = |x| $, \f \tx = (\tx_1, \tx_2, \cdots, \linebreak \tx_N)$ 
     $ = x/r$, \f \pa_j = \pa/\pa x_j \A and \f \gr = (\pa/\pa x_1, 
     \pa/\pa x_2, \cdots, \pa/\pa x_N)$. Then we set
\SP

        (1) \f k = k(x) = k(x, z) = [z\mu_0(x)]^{1/2}$, where the branch is 
     taken so that \f {\rm Im}\,k(x, z) \ge 0 $;       
\SP

        (2) \f a = a(x) = a(x, z) = {\rm Re}\,k(x, z)$;
\SP

        (3) \f b = b(x) = b(x, z) = {\rm Im}\,k(x, z)$;                           
\SP

        (4) \f \CD u = \gr u + \{(N-1)/(2r)\}\tx u - ik(x)\tx u$;                 
\SP

        (5) \f \CD_r u = \CD u \cdot \tx = \pa u/\pa r 
                                + \{(N-1)/(2r)\} u - ik(x) u$;                 
\BP

        Let \f u \in H^2(\RN)_{\rm{loc}} $. Then the restrictions
     \f u|_G \A and \f \pa_ju|_G $, \f j = 1, 2, \cdots, N $, of \f u \A
     and \f \pa_ju = \pa u/\pa x_j \A onto a smooth surface \f G \A 
     are defined as the traces of \f u \A and \f \pa_ju \A on \f G $, 
     respectively. Thus \f u|_G \A and \f \pa_ju|_G \A are considered to 
     belong to \f L_2(G)_{{\rm loc}} $.

        Let \f z \in \bC \A and let \f u \in H^2(\RN)_{{\rm loc}} $. Define 
     \f f \A by
$$
             f = -\,\mu_0(x)^{-1}\Delta u -\, zu
                   = \mu_0(x)^{-1}\big( -\,\Delta u - k^2u \big)   \T (3.2)
$$
     with \f k \A given by (1) of Notation 3.1. Now we are going to show 
     an identity which is an extension of Proposition 3.3 of [\JS] and will 
     be used throughout this section.
\BP

        {\bf Proposition 3.2.} {\it Let \f u \in H^2(\RN)_{\mb{loc}} \A 
     and let \f f \A be given by \r{\em (3.2)}$. Let \f \xi \A be a 
     real-valued, continuous function on \f [0, \infty) \A such that 
     \f \xi \A has piecewise continuous derivative. Set 
     \f \varphi(x) = \alpha(x)\xi(|x|)$, where \f \alpha \A is a simple 
     function which is constant on each 
     \f \Omega_{\ell} $. For \f 0 < r < R < \infty $, set
$$
         B_{rR} = \{ x \in \RN \ : \ r < |x| < R \ \},              \T (3.3)
$$
     Then we have     
$$
\split
    & {\dm \int_{B_{rR}} \big(b\varphi + \frac12\frac{\pa\varphi}{\pa r}\big)
                                                        |\CD u|^2\, dx 
          + \sum_{\ell\in L} \int_{\pa\Omega_{\ell}\cap B_{rR}}
       \varphi \mbox{\rm Im}\big\{\overline{k}\frac{\pa u}{\pa n}
                                 \overline{u}\big\}\, dS } \\  *[6pt]  
    & \hspace{4cm}  {\dm  + \int_{B_{rR}} \big(\frac{\varphi}r 
                                - \frac{\pa\varphi}{\pa r}\big)
                               (|\CD u|^2 - |\CD_ru|^2)\, dx } \\ *[6pt]
    & \hspace{4cm}  {\dm   + c_N \int_{B_{rR}} r^{-2}\big(\frac{\varphi}r 
              - 2^{-1}\frac{\pa\varphi}{\pa r} + b\varphi\big)|u|^2\, dx  } 
                                                                   \\ *[6pt]
    & {\dm = \mbox{\rm Re} 
                 \int_{B_{rR}} \varphi\mu_0(x)f\overline{\CD_ru}\, dx } 
                                                                  \\ *[6pt]
    &  \ \ \ \  {\dm  + 2^{-1}\sum_{\ell\in L} 
                         \int_{\pa\Omega_{\ell}\cap B_{rR}}
                     \varphi \big\{ \frac{(N-1)b}r + |k|^2 \big\}
                                      (\tx \cdot n)|u|^2\, dS } \\ *[6pt]
    &  \ \ \ \  {\dm  + 2^{-1} \int_{S_R} \varphi \big(2|\CD_ru|^2 - 
                         |\CD u|^2 - c_Nr^{-2}|u|^2 \big) \, dS } \\  *[6pt]
    &  \ \ \ \  {\dm - 2^{-1} \int_{S_r} \varphi(2|\CD_ru|^2 - |\CD u|^2 
                                               - c_Nr^{-2}|u|^2)\, dS, }
\endsplit                                                           \T (3.4)                                                             
$$
     where \f \Omega_{\ell} \A satisfies \r{\em (i), (ii)}$ of  
     \r{\em Assumption 2.1}$, \f \pa/\pa n \A in the integrand of the 
     surface integral over \f \pa\Omega_{\ell}\cap B_{rR} \A means the 
     directional derivative in the direction of the outward normal 
     \f n = n^{(\ell)} \A of \f \pa\Omega_{\ell}$, and}
$$     
              c_N = (N-1)(N-3)/4.                                \T (3.5)
$$  
\BP
      
        The proof will be omitted since it is essentially the same as 
     the proof of Proposition 3.3 of [\JS].    
\BP

        {\bf Theorem 3.3.} {\it Assume} \ Assumption 2.1. {\it Let 
     \f u \in H^2(\RN)_{{\rm loc}} \A be a solution of the homogeneous 
     equation 
$$
           -\, \mu_0(x)^{-1}\Delta u - \lambda u = 0 \qq (\lambda > 0)
                         			                                        		\T (3.6)
$$
     on \f \RN \A such that                           					
$$
     \liminf_{R\to\infty} 
         \int_{S_R} \big( \bigg|\frac{\pa u}{\pa r}\bigg|^2 
                        + |u|^2 \big)\, dS = 0,                     \T (3.7)
$$                        
     for \f N \ge 3 $, or
$$
     \liminf_{R\to\infty} 
         R^{\alpha}\int_{S_R} \big( \bigg|\frac{\pa u}{\pa r}\bigg|^2 
                        + |u|^2 \big)\, dS = 0                    \T (3.8)
$$                        
     with \f \alpha > 0 \A for \f N = 2 $, where
$$     
          S_R = \{ x \in \RN \ :\ |x| = R \}.                       \T (3.9)
$$
     Then \f u \A is identically zero.}
\BP

        Sketch of Proof. \ \ Theorme 3.3 can be proved by starting with 
     Proposition 3.2 and proceeding as in the proof of Theorem 3.2 
     of [\JS] (for \f N \ge 3 $) or proof of Theorem 7.1 of [\JS] 
     (for \f N = 2 $). Only difference here is that, instead of the last 
     inequality of (3.20) in [\JS], we have to use 
$$
     \sum_{\ell\in L} \int_{\pa\Omega_{\ell}\cap B_{rR}}
                     \varphi |k|^2(\tx \cdot n)|u|^2\, dS \le 0,    \T (3.10)
$$ 
     where \f n \A in the integrand is the unit outward normal 
     \f n^{(\ell)}(x) \A of \f \Omega_{\ell} \A at \f x $. In fact, it 
     follows from (i) and (ii) of Assumption 2.1 that
$$
\split
      & {\dm \sum_{\ell\in L} \int_{\pa\Omega_{\ell}\cap B_R} 
                     \varphi |k|^2(\tx \cdot n)|u|^2\, dS } \\ *[6pt]
      & \hspace{1cm} {\dm = \sum_{\ell\in L} \bigg( \int_{S_{\ell}^{(-)}\cap
                                                                      B_R} 
                        + \int_{S_{\ell}^{(+)}\cap B_R} \bigg)
                     \varphi |k|^2(\tx \cdot n)|u|^2\, dS } \\ *[6pt]
      & \hspace{1cm} {\dm = \lambda \sum_{\ell\in L} \int_{S_{\ell}^{(-)}\cap B_R}
                         \varphi\,\big(\nu_{\ell} - \nu_{\ell+1}\big)  
                            \,(\tx \cdot n^{(\ell)})|u|^2\, dS, }
\endsplit                                                         \T (3.11)
$$
    where we should note that we are dealing with a finite sum because of
    (i) of Assumption 2.1. Then (3.10) is obtained from (iii) of 
    Assumption 2.1.  \ \ \Box
\BP

        The following corollary guarantees the uniqueness of the inhomogeneous 
     equation
$$
         - \, \mu_0(x)^{-1}\Delta u - \lambda u = f               \T (3.12)
$$
     with one of the conditions
$$
         \M \CD_r^{(\pm)} u \M_{\delta-1, \, E_1} < \infty,         \T (3.13)
$$
     where \f \delta > 1/2 $,
$$
       \CD_r^{(\pm)} u = \pa u/\pa r + \{(N-1)/(2r)\} u \mp ik(x) u,
                                                                  \T (3.14)
$$     
$$
               E_R = \{ x \in \RN \ : \ |x| > R \},               \T (3.15)
$$
     and, for a measurable set \f G \A in \f \RN $,
$$
        \M v \M_{\delta-1, \, G}^2 
                = \int_{G} (1 + |x|)^{2(\delta-1)} |v(x)|^2 \, dx. \T (3.16) 
$$
\BP

        {\bf Corollary 3.4.} {\it Let \f \lambda > 0 \A and let 
     \f f \in L_2(\RN)_{{\rm  loc}} $. Then the solution 
     \f u \in H^2(\RN)_{{\rm  loc}} \A \linebreak of the equation} \ 
     (3.12) {\it with one of the radiation conditions in} \ (3.13) {\it is 
     unique.}
\BP

         The proof is the same as the proof of Corollary 3.8 of [\JS].

         Let \f L_{2, t}(\RN) \A be the weighted Hilbert space defined by 
     (1.8). Let the resolvent \f (H_0 - z)^{-1} \A of the operator \f H_0 \A
     be denoted by \f R_0(z) $. Now consider \f u \in X_0 \A defined by    
$$
\left\{ \split
          u = R_0(z)f,  \\ *[4pt]
          z = \lambda + i\eta  \qq (\lambda \ge 0, \eta \ne 0),  \\ *[4pt]
          f \in L_{2, \delta}(\RN).
\endsplit \right.                                            \T (3.17)
$$ 
     For \f 0 < c < d < \infty \A a subset
     \f J_{\pm}(c, \, d) \A of \f \bC \A are defined by
$$
\left\{ \split
      {\dm J_+(c, \, d) = \{ \ z = \lambda + i\eta \ : \ c \le \lambda \le d,
                                     \ 0 < \eta \le 1 \ \}, } \\  *[4pt]
      {\dm J_-(c, \, d) = \{ \ z = \lambda + i\eta \ : \ c \le \lambda \le d,
                                     \ -1 \le \eta < 0 \ \}.  }                              
\endsplit \right.                                                  \T (3.18)
$$
     
        In the next theorem we are going to evaluate the radiation 
     condition terms \f \CD u $. Here and in the sequel
     we agree that \f C = C(A, B, \cdots) \A in an inequality means a
     positive constant depending on \f A, B, \cdots $. Now we are evaluating
     the radiation condition term \f \CD u $.
\BP

        {\bf Theorem 3.5.} {\it Suppose that} \ Assumption 2.1 {\it holds.
     Let \f 1/2 < \delta \le 1 $. Let \f u \A be given 
     by} \ (3.17). 

        (i) \ {\it Let \f N \ge 3 $. Then there exists a constant 
     \f C = C(\delta, m_0, M_0) > 0 \A such that
$$
        \M \CD u \M_{\delta-1} \le C \M f \M_{\delta},
                                                                \T (3.19)   
$$
     where \f \CD u \A is as in} \ Notation 3.1, 
     {\it \f \M \ \M_t \A is the norm of \f L_{2, t}(\RN) $, and the constant 
     \f C(\delta) \A is independent of \f f \A and \f z \A satisfying}  
     \ (3.17).

        (ii) \ {\it Let \f N = 2 $. Let \f 0 < c < d < \infty \A and let 
     \f J_{\pm}(c, d) \A be as in {\rm (3.18)}. Let \f u \A be given by} \
     (3.17) {\it with \f z \in J_{+}(c, d) \cup J_{-}(c, d) $. Then there
     exists a positive constant \f C = C(\delta, c, d, m_0, M_0) \A such that
$$
        \M \CD u \M_{\delta-1,*} 
                  \le C \big( \M f \M_{\delta} + \M u \M_{-\delta} \big)
                                                                 \T (3.20)
$$
     where }
$$
       \M v \M_{t, *}^2 = \int_{B_1} |x||v(x)|^2 \, dx
                    + \int_{E_1} (1 + |x|)^{2t}|v(x)|^2 \, dx.   \T (3.21)
$$ 
\BP
       
        Sketch of Proof. \ \ We have only to proceed as in the proof of
     Theorems 4.1 and 7.2 of [\JS].  Set in (3.4) 
     \f \alpha(x) = 1/\sqrt{\mu_0} $, 
$$
     \xi(r) =
\left\{ \split
         r           \qq \q \hspace{1.5cm} (0 \le r \le 1), \\ *[4pt]
         2^{-(2\delta-1)}(1 + r)^{2\delta-1} \q (r \ge 1)
\endsplit \right.                                             \T (3.22)                  
$$  
     for \f N \ge 3 $, and
$$
      \xi(r) =
\left\{ \split
           {\dm \frac12 r^2 \qq \hspace{1.9cm} (r \le 1/2),} \\ *[4pt]
           {\dm \frac{1}{2^{2\delta}}(1 + r)^{2\delta-1}   \qq (r \ge 1) }.
\endsplit \right.                                                  \T (3.23)
$$
     for \f N = 2 $. Let the second term of the left-hand side of (3.4) be
     denoted by \f I_{L2} $. Then it is easy to see that
$$
     I_{L2} = \sum_{\ell\in L} \int_{\pa\Omega_{\ell}\cap B_{rR}}
       \varphi {\rm Im}\big\{\overline{k}\frac{\pa u}{\pa n}
                                 \overline{u}\big\}\, dS = 0   \T (3.24)
$$
     (cf. (3.11)). Similarly we see that the second term of the right-hand
     side of (3.4) is nonpositive. All other terms of (3.4) can be evaluated
     exactly in the same manner as in the proof of Theorems 4.1 and 7.2 of
     [\JS], which completes the proof.  \ \ \Box 
\BP

        Now that we have established the uniqueness of the solution of the 
     equation (3.12) with the radiation condition (Corollary 3.4) and the
     estimate of the radiation condition term (Theorem 3.5), we can show the 
     limiting absorption principle for \f H_0 \A by proceeding as in \S5,
     \S6, and \S7 of [\JS]. Let \f t \in \Ro $. The weighted Sobolev spaces 
     \f H_t^j(\RN) $, \f j = 1, 2 $, are defined as the completion of 
     \f \Con(\RN) \A by the norms
$$
        \M u \M_{1,t} 
          = \bigg[ \int_{\RN} (1 + r)^{2t} \big( |\gr u|^2 + 
                              |u(x)|^2 \big) \, dx \bigg]^{1/2},   \T (3.25)
$$
     and
$$
        \M u \M_{2,t} 
          = \bigg[ \int_{\RN} (1 + r)^{2t} \sum_{|\gamma|\le 2}
                        |\pa^{\gamma}u|^2 \, dx \bigg]^{1/2},      \T (3.26)                           
$$
     respectively, where
$$
\left\{ \split
       {\dm \gamma = (\gamma_1, \gamma_2, \cdots, \gamma_N), } \\ *[4pt]
       {\dm |\gamma| = \gamma_1 + \gamma_2 + \cdots + \gamma_N, } \\ *[4pt]
       {\dm \pa^{\gamma}u = (\pa_1)^{\gamma_1} \cdots (\pa_N)^{\gamma_N}u
            \q (\pa_j = \pa/\pa x_j). }
\endsplit \right.                                                 \T (3.27)
$$
    The inner product and norm of \f H_t^j(\RN) \A will be denoted by
    \f (\ , \ )_{j,t} \A and \f \M \ \M_{j,t} $. For an operator \f T $,\, 
    the operator norm in \f \bB(H_s^j(\RN), \, H_t^{\ell}(\RN)) $ will
    be denoted by  \f \M T \M_{(j,s)}^{(\ell,t)} $, where \f j, \ell = 0, 1, 
    2 $, \f s, t \in \Ro $, and we set
$$
         H_s^0(\RN) = L_{2,s}(\RN).                               \T (3.28)
$$         
     Let \f D_{\pm} \subset \bC \A be given by
$$
\left\{ \split
  D_+ 
     = \{ \, z = \lambda + i\eta \, : \, \lambda > 0, \, \eta \ge 0 \, \}, \\ *[4pt]       
  D_-
     = \{ \, z = \lambda + i\eta \, : \, \lambda > 0, \, \eta \le 0 \, \}.     
\endsplit \right.                                               \T (3.29)
$$
     Also, for \f 0 < c < d < \infty $, let \f J_{\pm}(c, d) \A
     be as in (3.18). The closure  \f \overline{J}_{\pm}(c, d) \A 
     \linebreak are given by
$$
\left\{ \split
      \overline{J}_+(c, \, d) = \{ \ z = \lambda + i\eta \ : \ 
          c \le \lambda \le d, \ 0 \le \eta \le 1 \ \} \subset D_+,  
                                                              \\  *[4pt]   
      \overline{J}_-(c, \, d) = \{ \ z = \lambda + i\eta \ : \ 
          c \le \lambda \le d, \ - \, 1 \le \eta \le 0 \ \} \subset D_-.                                
\endsplit \right.                                                  \T (3.30)
$$ 
     For \f \lambda > 0$, let 
$$
            R_{0\pm}(\lambda) 
                   = \lim_{\eta\downarrow 0} R_0(\lambda \pm i\eta),   
                                                                  \T (3.31)
$$
     and extend the resolvent \f R_0(z) \A on \f D_{\pm} \A by
$$
     R_0(\lambda+i\eta) = 
\left\{ \split
    R_0(\lambda+i\eta) \qq (\lambda > 0, \, \eta > 0), \\  *[4pt]      
    R_{0+}(\lambda) \ \ \qq \ \ (\lambda > 0, \, \eta = 0)    
\endsplit \right.                                               \T (3.32)
$$
     for \f z \in D_+ \A and                   
$$
     R_0(\lambda+i\eta) = 
\left\{ \split
    R_0(\lambda+i\eta) \qq (\lambda > 0, \, \eta < 0), \\ *[4pt]       
    R_{0-}(\lambda) \ \ \qq \ \ (\lambda > 0, \, \eta = 0)    
\endsplit \right.                                               \T (3.33)
$$
     for \f z \in D_- $. Then we have
\BP

        {\bf Theorem 3.6.} {\it Suppose that} \ Assumption 2.1 {\it holds. 
     Let \f 1/2 < \delta \le 1 \A. }
\SP

        (i) \ {\it Then the limits} \ (3.31) {\it is well-defined in
     \f \bB(L_{2,\delta}(\RN), \, H_{-\delta}^2(\RN)) $, and the extended
     resolvent \f R_0(z) \A is a 
     \f \bB(L_{2,\delta}(\RN), \, H_{-\delta}^2(\RN))$-valued continuous 
     function on each of \f D_+ \A and \f D_- $.} 
\SP

        (ii) \ {\it For any \f z \in D_+ \A {\rm [}\,or \f D_- $\,{\rm ]}, 
     \f R_0(z) \A is a compact operator from \f L_{2,\delta}(\RN) \A 
     into \f H_{-\delta}^1(\RN) $.}     
\SP

        (iii) \ {\it The selfadjoint operator \f H_0 \A is absolutely 
     continuous on the interval \f (0, \infty) $. The operator \f H_0 \A has
     neither point spectrum nor singular continuous spectrum.}
\SP

        (iv) \ {\it  For \f 0 < c < d < \infty \A there exists a  
     constant \f C = C(c, d, \delta, m_0, M_0) > 0 \A such that, for 
     \f z \in \overline{J}_+(c, \, d) \cup \overline{J}_-(c, \, d) $,
$$   
\left\{ \split
      {\dm \int_{E_s}(1 + r)^{-2\delta} 
            \big( |\gr R_0(z)f|^2 + |k|^2|R_0(z)f|^2 \big) \, dx } \\  *[5pt]
      \hspace{3.5cm} {\dm  \le C^2 (1 + s)^{-(2\delta-1)} 
                                          \M f \M_{\delta}^2 }\\ *[4pt]
       \ \ \ \ \ \ \ \ \ \ \ \ \ \ \ \ \ \ \ \ \ \ \ \ \ \ \ \ \ \ \ \ \
         \ \ \ \ \ \ \ \ {\dm (s \ge 1, \ f \in L_{2,\delta}(\RN)), } 
                                                                \\  *[4pt] 
     {\dm \M \CD R_0(z)f \M_{\delta-1} \le C \M f \M_{\delta} \q
             ( f \in L_{2,\delta}(\RN)), }     
\endsplit \right.                                              \T (3.34)
$$
     where, for \f \lambda \in D_+ \cap (0, \, \infty) \A 
     or \f D_- \cap (0, \, \infty)$, \f \CD u \A should 
     be interpreted as}
$$
            \CD u = 
\left\{ \split
         {\dm  \CD^{(+)} u = \gr u + \{(N-1)/(2r)\}\tx u - ik(x)\tx u \q
                         (\lambda \in D_+ \cap (0, \, \infty)), } \\ *[5pt]
         {\dm  \CD^{(-)} u = \gr u + \{(N-1)/(2r)\}\tx u + ik(x)\tx u \q
                          (\lambda \in D_- \cap (0, \, \infty)).} \\  *[5pt]
\endsplit \right.                                               \T (3.35)
$$ 
\SP

        (v) \ {\it Let \f N \ge 3 $. Then there exists a constant 
     \f C = C(\delta, m_0, M_0) > 0 \A such that }
$$
\left\{ \split
     {\dm \int_{E_s} (1 + r)^{-2\delta} |R_0(z)f|^2 \, dx 
         \le \frac{C^2}{|z|}(1 + s)^{-2(2\delta-1)} \M f \M_{\delta}^2 } 
                                                             \\ *[5pt]
       \ \ \ \ \ \ \ \ \ \ \ \ \ \ \ \ \ \ \ \ \ \ \ \ \ \ \ \ \ \ \ \
           \ \ \ \ \ \ \ \ \ \ \ \ \ \  (s \ge 0, \ f \in L_{2,\delta}),  
                                                              \\ *[5pt] 
     {\dm \M R_0(z) \M_{(0,\delta)}^{(0,-\delta)} \le \frac{C}{\sqrt{|z|}}
                                    \q (z \in D_+ \cup D_-), } \\ *[5pt]
     {\dm \M \CD R_0(z)f \M_{\delta-1} \le C \M f \M_{\delta} \q
             (z \in D_+ \cup D_-, \ f \in L_{2,\delta}). }                            
\endsplit \right.                                              \T (3.36)
$$    
\BP

        Finally we are going to prove a modification of Theorem 3.5, where
     the range of \f \delta \A is slightly wider. This modification will be
     useful in the next section.
\BP

        {\bf Proposition 3.7.} {\it Let} \ Assumption 2.1 {\it be satisfied. 
     Let \f 1/2 < \delta < 3/2 $. Let \f f \in L_{2, \delta}(\RN) \A and
     let \f z \in D_+ \cup D_- $. Let \f u = R_0(z)f $.} 

        (i) \ {\it Let \f N \ge 3 $. Then there exists a constant 
     \f C = C(\delta, m_0, M_0) > 0 \A such that
$$
        \M \CD_r u \M_{\delta-1} \le C \M f \M_{\delta} \q
                            (f \in L_{2, \delta}(\RN), \ z \in D_+ \cup D_-)
                                                                \T (3.37)   
$$
     where \f \CD_r u \A is given in} \ Notation 3.1, (6), {\it and
     for \f \lambda \in D_+ \cap (0, \, \infty) \A 
     [or \f D_- \cap (0, \, \infty)$], \f \CD_r u \A should be interpreted 
     as \f \CD_r^{(+)} \A [\,or \f \CD_r^{(-)} \,$].} 
\SP

        (ii) \ {\it Let \f N = 2 $. Let \f 0 < c < d < \infty $. Then 
     there exists \f C = C(\delta, c, d,  m_0, M_0) > 0 \A such that
$$
        \M \CD_r u \M_{\delta-1,*} \le C \M f \M_{\delta} \q
                            (f \in L_{2, \delta}(\RN), \ 
               z \in \overline{J}_{+}(c, d) \cup \overline{J}_{-}(c, d)),
                                                                \T (3.38)   
$$    
     where \f \M \CD_r u \M_{\delta-1,*} \A is given by} \ (3.21).
\BP

        Proof. \ \ In view of the continuity of the extended resolvent
     of \f R_0(z) $, we only have to prove (3.37) and (3.38) for non real
     \f z $. Then we should note that we have 
     \f  u = R_0(z)f \in H_{\delta}^2(\RN) \A (cf., e.g., [\Sc], Lemma 2.1).
     As in the proof of Theorem 3.5, we start with (3.4) with \f \xi \A 
     given by (3.22) or (3.23) and \f \alpha(x) = 1/\sqrt{\mu_0} $. Let the 
     {\it j}-th term of the left-hand side be denoted by \f I_{Lj} $, where 
     \f j = 1, 2, 3, 4 $. Then we have \f  I_{L2} = 0 $, and, since 
$$
     \frac{\varphi}r - 2^{-1}\frac{\pa\varphi}{\pa r} \ge 0    \T (3.39)
$$
     and \f |\CD u| \ge |\CD_r u| $, we have
$$
\split
   {\dm I_{L1} + I_{L3} \ge \int_{B_{rR}} \frac12\frac{\pa\varphi}{\pa r}
                                                        |\CD u|^2\, dx 
                            + \int_{B_{rR}} \big(\frac{\varphi}r 
                                - \frac{\pa\varphi}{\pa r}\big)
                               (|\CD u|^2 - |\CD_ru|^2)\, dx } \\ *[6pt]
    \hspace{1.7cm} {\dm = \int_{B_{rR}} \big(\frac{\varphi}r 
                    - \frac12 \frac{\pa\varphi}{\pa r}\big) |\CD u|^2\, dx 
                    - \int_{B_{rR}} \big(\frac{\varphi}r 
                                - \frac{\pa\varphi}{\pa r}\big)
                                           |\CD_ru|^2\, dx } \\ *[6pt]
    \hspace{1.7cm} {\dm \ge \int_{B_{rR}} \big(\frac{\varphi}r 
                    - \frac12 \frac{\pa\varphi}{\pa r}\big) |\CD_r u|^2\, dx 
                    - \int_{B_{rR}} \big(\frac{\varphi}r 
                                - \frac{\pa\varphi}{\pa r}\big)
                                           |\CD_ru|^2\, dx } \\ *[6pt]
      \hspace{1.7cm} {\dm = \int_{B_{rR}}  
                    \frac12 \frac{\pa\varphi}{\pa r} |\CD_r u|^2\, dx.}
\endsplit                                                         \T (3.40)
$$ 
     As for the fourth term \f I_{L4} $, we have \f I_{L4} \ge 0 \A for 
     \f N \ge 3 $, and for \f N = 2 \A we have, as in (7.19) and (7.20) of 
     [\JS], 
$$
\split
     {\dm -\,I_{L4} \le C_1 \M u \M_{\delta-2}^2 
                          + C_2 \int_{\Rb} |\eta||u|^2 \, dx } \\ *[6pt]
      \hspace{1.1cm} {\dm \le C_1 \M u \M_{\delta-2}^2 
                               + C_3(|f|, |u|)_0 }  \\ *[6pt]
      \hspace{1.1cm} {\dm \le C_1 \M u \M_{\delta-2}^2 
                               + \frac{C_3}2 (\M f \M_{\delta}^2
                                     + \M u \M_{-\delta}^2) }   
\endsplit                                                        \T (3.41)
$$
     with \f C_1 = C_1(\delta, m_0) $, \f C_2 =C_2(c, d) $, and
     \f C_3 = C_3(\delta, c, d, m_0, M_0) $. Here we can evaluate
     \f \M u \M_{\delta-2} \A as
$$
\left\{ \split
       {\dm \delta \in (1/2, 1] \Longrightarrow \M u \M_{\delta-2}
                  \le \M u \M_{-\delta} \le C'\M f \M_{\delta}, } \\ *[4pt]
       {\dm \delta \in (1, 3/2) \Longrightarrow \M u \M_{\delta-2}
                     \le C''\M f \M_{2-\delta} \le C''\M f \M_{\delta}, } 
\endsplit \right.                                                 \T (3.42)
$$
     with \f C' = C'(\delta, c, d, \mu_0) \A and 
     \f C'' = C''(\delta, c, d, \mu_0) $. We can evaluate the right-hand 
     of (3.4) by proceeding as in the proof of Theorems 4.1 and 7.2 of [\JS]
     Therefore, letting \f r \to 0 \A and \f R \to \infty $, where we have 
     to use the fact that \f  u = R_0(z)f \in H_{\delta}^2(\RN) $, 
     we obtain (3.37) and (3.38), respectively. This completes the proof.
     \ \ \Box 
\BP
        
$$
\ \ \ 
$$

        {\bf 4. \ The point spectrum for \f H $}              
\BP

        Throughout this and the following sections we shall assume that 
     Assumptions 2.1 and 2.2 are satisfied. Let the operator \f H \A be as 
     in \S2. Let \f \sigma_p(H) \A be the set of all eigenvalues of \f H $, 
     and let \f V_p(H) \A be the set of all eigenvectors of \f H $, i.e.,     
$$
     V_p(H) = \{ u \in H^2(\RN) : u \ne 0, \, (H - \lambda)u = 0
                         \ \ {\rm with} \ \ \lambda \in \sigma_p(H) \}.
                                                               \T (4.1)
$$
     Now we need to introduce a set \f \widetilde{\sigma}_p^{(\pm)}(H) \A of
     the extended eigenvalues of \f H \A and a set 
     \f \widetilde{V}_p^{(\pm)}(H) \A of the extended eigenvectors of \f H $.
\BP

       {\sc Definition 4.1.} Let \f 1/2 < \delta < 1/2 + \epsilon 
     \A if \f \mu_1 \A is short-range, i.e., \f \mu_1 \A satisfies (2.13), 
     and let \f 1/2 < \delta < (1 + \epsilon)/2 \A if \f \mu_1 \A is 
     long-range, i.e., \f \mu_1 \A satisfies (2.14), where \f \epsilon \A 
     is given in Assumption 2.2. Denote by 
     \f \widetilde{\sigma}_p^{(+)}(H) \A [or 
     \f \widetilde{\sigma}_p^{(-)}(H) $] the set of all \f \lambda > 0 \A 
     such that there exists a function \f u \A satisfying
$$
\split
     {\rm (i)}  \ \ \ u \in H^2(\RN)_{{\rm loc}}, \ u \ne 0, \\
     {\rm (ii)} \ \  u \in L_{2, -\delta}(\RN), \\
     {\rm (iii)} \ \ \M \CD^{(+)}u \M_{\delta-1, E_1} < \infty, 
          \ {\rm [or} \ \M \CD^{(-)}u \M_{\delta-1, E_1} < \infty {\rm ]}, \\
     {\rm (iv)} \ \ \ -\,\mu(x)^{-1}\Delta u - \lambda u = 0,
\endsplit                                                     \T (4.2)
$$
     where \f \CD^{(\pm)} u \A are given by (3.35),
     and \f k = k(x) = [\lambda\mu_0(x)]^{1/2} \A is as in \S3. Let 
     \f \widetilde{V}_p^{(\pm)}(H) \A be the set of all \f u \in X \A which 
     satisfy (4.2) with \f \lambda \in \widetilde{\sigma}_p^{(\pm)}(H) $. 
     We call \f u \in \widetilde{V}_p^{(\pm)}(H) \A which satisfies (4.2) 
     an extended eigenvector of \f H \A associated with the extended 
     eigenvalue \f \lambda $.
\BP

        Since \f 0 \not\in \sigma_p(H) $, we have 
     \f V_p(H) \subset \widetilde{V}_p^{(\pm)}(H) \A and
     \f \sigma_p(H) \subset \widetilde{\sigma}_p^{(\pm)}(H) \A by definition. 
     In this section we are going to prove that 
$$
    \sigma_p(H) = \widetilde{\sigma}_p^{(+)}(H) 
                   = \widetilde{\sigma}_p^{(-)}(H),             \T (4.3)
$$
     and \f \sigma_p(H) \A is a discrete set on \f (0, \infty) $.
\BP

        {\bf Proposition 4.2.} {\it Assume that} \ Assumptions 2.1 {\it and}
     \ 2.2 {\it hold. Let \f u \in H^2(\RN)_{{\rm loc}} \A be 
     a solution of the equation \f -\,\mu(x)^{-1}\Delta u - \lambda u = 0 $.
     Let \f \varphi(x) = \xi(|x|) \A and let \f \xi(r) \A satisfy the 
     following} \ (a) $\sim$ (c)':

        (a) \ {\it  \f \xi \A is a nonnegative, continuous function on 
     \f (0, \infty) $,}

        (b) \ {\it  \f \xi \A has a piecewise continuous derivative 
     \f \xi' \A such that
$$
         \xi'(r) \ge 0,                                             \T (4.4)
$$
     and
$$         
         {\dm \frac{\xi(r)}r - \frac12\xi'(r) \ge 0 }               \T (4.5)
$$
     for almost all \f r > 0 $.} 

        (c) \ {\it If \f N \ge 3 $,}
$$
            \xi(r) = O(r) \qq (r \to 0).                            \T (4.6)
$$
        
        (c)' \ {\it If \f N = 2 $,}
$$
\left\{ \split
         \xi(r) = O(r^2) \qq (r \to 0),  \\ *[4pt]
         \xi'(r) = O(r) \qq  (r \to 0).
\endsplit \right.                                                  \T (4.7)
$$

        (i) \ {\it Suppose that \f \mu_1 \A is short-range in the sense of} 
     \ (2.13). {\it Then there exists a constant 
     \f C = C(\lambda, \mu_0) > 0 \A such that, for \f R > 1 $,
$$
\split
    {\dm \int_{B_R} \frac{\pa\varphi}{\pa r}|u|^2 \,dx
            \le 2m_0^{-1} \int_{B_R} \varphi|\mu_1||u||\CD_r^{(\pm)}u|\,dx } 
                                                                \\ *[5pt]
    \hspace{3cm} {\dm + \ C\xi(R) \int_{S_R} |\CD_r^{(\pm)}u|^2 \,dS }
\endsplit                                                         \T (4.8)
$$
     for \f N \ge 3 $, or
$$
\split
    {\dm \int_{B_R} \frac{\pa\varphi}{\pa r}|u|^2 \,dx
          \le 2m_0^{-1} \int_{B_R} \varphi|\mu_1||u||\CD_r^{(\pm)}u|\,dx } 
                                                                  \\  *[5pt]
     \hspace{3cm} {\dm + \,(2m_0\lambda)^{-1} \int_{B_R} r^{-2}
          \bigg( \frac{\varphi}r - 2^{-1}\frac{\pa\varphi}{\pa r} \bigg)
                                                       |u|^2\,dx } \\ *[5pt]
    \hspace{3cm} {\dm + \ C\xi(R) \int_{S_R} |\CD_r^{(\pm)}u|^2 \,dS }
\endsplit                                                           \T (4.9)
$$
     for \f N = 2 $.}

        (ii) {\it  Suppose that \f \mu_1 \A is long-range in the sense of}
     \ (2.14). {\it Then the relation} \ (4.8) {\it or} \ (4.9) {\it holds 
     again with the first term \f K_{R1} \A of the right-hand side of } \
     (4.8) {\it or} \ (4.9) {\it replaced by}
$$
    K_{R1}' = m_0^{-1}\, \int_{B_R} \bigg( \bigg|\frac{\pa\varphi}{\pa r}
                  \mu_1\bigg| + \bigg|\frac{\pa\mu_1}{\pa r}\varphi\bigg| 
                  \bigg) |u|^2 \, dx.
                                                                 \T (4.10)
$$
\BP 

        Proof. \ \ (I) We shall prove (4.8) and (4.9) for 
     \f \CD_r^{(+)}u $. These formulas for \f \CD_r^{(-)}u \A can be 
     proved in quite a similar way. For the sake of simplicity of notation 
     we set \f \CD_r^{(+)}u = \CD_ru \A and \f \CD^{(+)}u = \CD u $. 
     Since we have from (iv) of (4.2)
     \f -\,\Delta u - \lambda\mu_0(x)u = \lambda\mu_1(x)u $,          
     we can apply the formula (3.4) with \f f \A and \f z \A replaced by 
     \f \lambda\mu_0(x)^{-1}\mu_1(x)u \A and \f \lambda $, respectively, 
     to obtain, for \f 0 < r < 1 < R < \infty $,
$$
\split
    & {\dm \int_{B_{rR}} \frac12\frac{\pa\varphi}{\pa r}|\CD u|^2\, dx 
          + \sum_{\ell\in L} \int_{\pa\Omega_{\ell}\cap B_{rR}}
       \varphi {\rm Im}\big\{\overline{k}\frac{\pa u}{\pa n}
                                 \overline{u}\big\}\, dS } \\  *[6pt]  
    & \hspace{2cm}  {\dm  + \int_{B_{rR}} \big(\frac{\varphi}r 
                                - \frac{\pa\varphi}{\pa r}\big)
                               (|\CD u|^2 - |\CD_ru|^2)\, dx } \\ *[6pt]
    & \hspace{2cm}  {\dm   + c_N \int_{B_{rR}} r^{-2}\big(\frac{\varphi}r 
              - 2^{-1}\frac{\pa\varphi}{\pa r}\big)|u|^2\, dx  } \\ *[6pt]
\endsplit                                                         \T (4.11)
$$
$$
\split
    & {\dm = {\rm Re} 
           \int_{B_{rR}} \lambda\varphi\mu_1(x)u\overline{\CD_ru}\, dx } 
                                                                  \\ *[6pt]
    &  \ \ \ \  {\dm  + \ 2^{-1}\sum_{\ell\in L} 
                      \int_{\pa\Omega_{\ell}\cap B_{rR}}
                     \varphi k^2 (\tx \cdot n)|u|^2\, dS } \\ *[6pt]
    &  \ \ \ \  {\dm  + \ 2^{-1} 
                    \int_{S_R} \varphi \big(2|\CD_ru|^2 - |\CD u|^2 
                                   - c_Nr^{-2}|u|^2 \big) \, dS } \\  *[6pt]
    &  \ \ \ \  {\dm - \ 2^{-1} \int_{S_r} \varphi(2|\CD_ru|^2 - |\CD u|^2 
                                               - c_Nr^{-2}|u|^2)\, dS. }
\endsplit                                                                                                                       
$$
\SP

        (II) Suppose that \f \mu_1 \A is short-range. Proceeding as in the 
     proof of [\JS], Theorems 3.2 or 7.1, we have
$$
\split    
     {\dm \int_{B_R} \frac12 \frac{\pa\varphi}{\pa r} |\CD u|^2\, dx  
                + \sum_{\ell\in L} \int_{\pa\Omega_{\ell}\cap B_R}
                    \varphi  \mbox{\rm Im}\big\{k\frac{\pa u}{\pa n}
                                 \overline{u}\big\}\, dS } \\ *[6pt]
     \hspace{1.5cm} 
         {\dm \ge \int_{B_R} \frac12 \frac{\pa\varphi}{\pa r} k^2 |u|^2\, dx 
              + \int_{B_R} \frac12 \frac{\pa\varphi}{\pa r} ( |\gr u|^2 
                             - \bigg|\frac{\pa u}{\pa r}\bigg|^2)\, dx  } 
                                                                  \\ *[8pt]
     \hspace{3cm} 
         {\dm - \ \xi(R) \int_{S_R} \mbox{\rm Im}\big\{ k\frac{\pa u}{\pa r}
                                 \overline{u} \big\}\, dS, }       
\endsplit                                                          \T (4.12)                                                 
$$
     where  we should note that all the integrals in (4.12) is well-defined 
     even in the case of \f N = 2 \A because of (4.7). Therefore it follows 
     from (4.11) and (4.12) that
$$
\split
      {\dm \int_{B_{R}} \frac12 k^2 \frac{\pa\varphi}{\pa r}|u|^2\, dx } 
                                                                   \\ *[6pt]
         \hspace{2cm}
           {\dm + \int_{B_{rR}} \big(\frac{\varphi}r 
            - 2^{-1}\frac{\pa\varphi}{\pa r}\big)
               (|\gr u|^2 - |\frac{\pa u}{\pa r}|^2)\, dx  } \\ *[6pt]
         \hspace{2cm}
           {\dm + \ c_N \int_{B_{rR}} r^{-2}\big( \frac{\varphi}r 
                  - 2^{-1}\frac{\pa\varphi}{\pa r}\big)|u|^2\, dx } \\ *[6pt]
\endsplit
$$
$$
\split                  
     \hspace{1cm}
       {\dm \le {\rm Re} 
      \int_{B_{rR}} \lambda\varphi\mu_1(x)u\overline{\CD_ru}\, dx } \\ *[6pt]
         \hspace{2cm}
           {\dm + \ 2^{-1}\sum_{\ell\in L} \int_{\pa\Omega_{\ell}\cap B_{rR}}
                     \varphi |k|^2(\tx \cdot n)|u|^2\, dS } \\ *[7pt]  
        \hspace{2cm} 
           {\dm + \ 2^{-1} \xi(R) \int_{S_R} (2|\CD_ru|^2 - |\CD u|^2 
                    - c_Nr^{-2}|u|^2)\, dS } \\ *[7pt]                
        \hspace{2cm} 
           {\dm - \ 2^{-1} \xi(r) \int_{S_r} (2|\CD_ru|^2 - |\CD u|^2 
                    - c_Nr^{-2}|u|^2)\, dS, } \\ *[7pt]   
        \hspace{2cm} 
           {\dm  + \int_{B_r} \frac12 \frac{\pa\varphi}{\pa r}|\CD u|^2\, dx 
          + \sum_{\ell\in L} \int_{\pa\Omega_{\ell}\cap B_r}
                \varphi{\rm Im}\big( k\frac{\pa u}{\pa n}
                                 \overline{u}\big) \, dS  } \\ *[7pt]                                                              
         \hspace{2cm} 
           {\dm + \ \xi(R) \int_{S_R} \mbox{\rm Im}\big( k\frac{\pa u}{\pa r}
                                 \overline{u} \big) \, dS. }   
\endsplit                                                         \T (4.13)
$$
     It follows from (2.11) that the second term of the right-hand side of
     (4.13) is nonpositive. Therefore we can drop it from the right-hand 
     side of (4.13). Further, we see from (4.5) in (b) that the second term 
     of the left-hand side of (4.13) is nonnegative, and it can be dropped, 
     too. Then, letting \f r \to 0 \A along an appropriate sequence
     \f \{ r_n \} $, we obtain 
$$
\split
      {\dm \lambda \int_{B_{R}} \frac12 \mu_0(x) \frac{\pa\varphi}{\pa r}
                                                                 |u|^2\, dx 
                   + c_N \int_{B_R} r^{-2}\big( \frac{\varphi}r 
                  - 2^{-1}\frac{\pa\varphi}{\pa r}\big)|u|^2\, dx } \\ *[6pt]
           \hspace{1cm}
           {\dm \le 
              \lambda \int_{B_R} \varphi|\mu_1(x)||u||\CD_ru|\, dx } 
                                                                   \\ *[6pt]
           \hspace{2cm} 
            {\dm + \ 2^{-1} \xi(R) \int_{S_R} (2|\CD_ru|^2 - |\CD u|^2 
                    - c_Nr^{-2}|u|^2)\, dS   }  \\ *[7pt]             
        \hspace{2cm} 
            {\dm + \ \xi(R) \int_{S_R} {\rm Im}\big( k\frac{\pa u}{\pa r}
                                 \overline{u} \big) \, dS. }  \\ *[7pt]
        \hspace{1cm}
           {\dm \le 
              \lambda \int_{B_{R}} \varphi|\mu_1(x)||u||\CD_ru|\, dx } 
                                                                   \\ *[7pt]
        \hspace{2cm}
           {\dm + \ C\xi(R) \int_{S_R} \bigg( 
                         \bigg|\frac{\pa u}{\pa r}\bigg|^2
                               + |u|^2 \bigg) \, dS. } 
\endsplit                                                    \T (4.14)         
$$
     Here we noticed from (c) and (c)' that 
$$
\left\{ \split
     {\dm \liminf_{r\to 0} \ \xi(r) \int_{S_r} (2|\CD_ru|^2 - |\CD u|^2 
                    - c_Nr^{-2}|u|^2)\, dS = 0, } \\ *[7pt] 
     {\dm \ 2^{-1} \xi(R) \int_{S_R} (2|\CD_ru|^2 - |\CD u|^2 
                    - c_Nr^{-2}|u|^2)\, dS   }  \\ *[7pt]             
        \hspace{2cm} 
            {\dm + \ \xi(R) \int_{S_R} {\rm Im}\big( k\frac{\pa u}{\pa r}
                                 \overline{u} \big) \, dS. }  \\ *[7pt]
       \hspace{2cm} {\dm \le  C\xi(R) \int_{S_R} \bigg( 
                         \bigg|\frac{\pa u}{\pa r}\bigg|^2
                               + |u|^2 \bigg) \, dS } 
\endsplit \right.                                              \T (4.15)
$$
     with a constant \f C $. Since \f c_N \ge 0 \A for \f N \ge 3 \A and 
     \f c_N = -1/4 \A for  \f N = 2 $, (4.8) and (4.9) follows from (4.14) 
     if we can prove
$$
     \int_{S_R} \bigg( \bigg|\frac{\pa u}{\pa r}\bigg|^2
                               + |u|^2 \bigg) \, dS
            \le C'\int_{S_R} |\CD_r^{(\pm)}u|^2 \,dS             \T (4.16)
$$
     with a constant \f C' $.
\SP

        (III) (Proof of (4.16).) Multiply both sides of the equation
     \f -\,\Delta u - \mu(x)\lambda u = 0 \A by \f \overline{u} $, 
     integrating over \f B_R $, use partial integration and take the 
     imaginary part to obtain
$$
         \int_{S_R} {\rm Im}\big\{ \frac{\pa u}{\pa r}
                                 \overline{u} \big\}\, dS = 0.      \T (4.17)
$$
     Now we can proceed as in the proof of [\JS], Theorem 3.7 to obtain 
     (4.16).
\SP

        (IV) Suppose that \f \mu_1 \A is long-range. We have only to evaluate
     the first term \f K_{R1}'' \A of the right-hand side of (4.11). In fact, we have
     by partial integration
$$
\split
      {\dm K_{R1}'' = \lambda \, {\rm Re} \, \int_{B_{rR}} \varphi\mu_1(x)u
                                       \overline{\CD_ru}\, dx } \\   *[6pt]
      \hspace{.9cm}{\dm  = \frac{\lambda}2 \, \int_{B_{rR}} \varphi\mu_1(x) 
                   \big( \frac{\pa |u|^2}{\pa r} + (N-1)r^{-1}|u|^2 \big)
                                                         \, dx } \\   *[6pt]
      \hspace{.9cm}{\dm = - \frac{\lambda}2 \, 
                  \int_{B_{rR}} \frac{\pa (\varphi\mu_1)}{\pa r}|u|^2 \, dx
                    + \frac{\lambda}2 \, 
                    \int_{S_{R}} \varphi\mu_1|u|^2 \, dS }  \\   *[6pt]
       \hspace{1.9cm} {\dm - \frac{\lambda}2 \, 
                    \int_{S_{r}} \varphi\mu_1|u|^2 \, dS } 
\endsplit                                                          \T (4.18)
$$
     Since the third term of the right-hand side of (4.18) converges to 
     \f 0 \A as \f r \to 0 $, we see that \f K_{R1}' \A in (4.10) can replace the 
     first term \f K_{R1} \A of the right-hand side of (4.8) or (4.9), which 
     completes the proof. \ \ \Box
\BP

        {\bf Proposition 4.3.} {\it  Assume that} \ Assumptions 2.1 
     {\it and} \ 2.2 {\it hold. Suppose that \f \mu_1 \A is short-range. 
     Suppose that \f u \in \widetilde{V}_p^{(+)}(H) \A [\,or \f u \in 
     \widetilde{V}_p^{(-)}(H) \,$] with an extended eigenvalue 
     \f \lambda \in \widetilde{\sigma}_p^{(+)}(H) \A [\,or 
     \f \lambda \in \widetilde{\sigma}_p^{(-)}(H) \, $] such that
$$
\left\{ \split
              u \in L_{2, -\delta+j\epsilon}(\RN),  \\ *[5pt]    
              - \delta + j\epsilon  \le 0   
\endsplit \right.                                                 \T (4.19)
$$
     with a nonnegative integer \f j $. Then we have 
$$
\left\{ \split
   u \in L_{2, -\delta+(j+1)\epsilon}(\RN), \\ *[4pt]
   \CD_r^{(+)}u \in L_{2, -\delta+(j+1)\epsilon}(\RN) \ \
             [ \, {\rm or} \ \CD_r^{(-)}u \in 
                        L_{2, -\delta+(j+1)\epsilon}(\RN) \, ],
\endsplit \right.                                                 \T (4.20)
$$
     and
$$
\left\{ \split
       \M u \M_{-\delta+(j+1)\epsilon} 
             \le \sqrt{\lambda}c_1C_j^{(N)}\M u \M_{-\delta+j\epsilon} 
                                             \q (N \ge 3), \\ *[4pt]
       \M u \M_{-\delta+(j+1)\epsilon} 
             \le c_1C_j^{(2)}\M u \M_{-\delta+j\epsilon} 
                                             \q (N = 2),             
\endsplit \right.                                                 \T (4.21)
$$
     where \f C_j^{(N)} = C_j^{(N)}(\delta, \epsilon, m_0, M_0) \A for 
     \f N \ge 3 \A and  
     \f C_j^{(2)} = C_j^{(2)}(\lambda, \delta, \epsilon, m_0, M_0) $.}
\BP

        Proof. \ \ (I) We shall prove (4.20) for \f \CD_r^{(+)}u $, and set 
     \f \CD_r^{(+)}u = \CD_ru $. Then \f u \A satisfies the equation 
     \f (- \Delta - \lambda\mu_0)u =\lambda\mu_1u \A and 
     the radiation condition \f \M \CD_r^{(+)}u \M_{\delta-1,E_1} < \infty $,
     i.e.,
$$
            u = \lambda R_0(\lambda)\big( \mu_0^{-1}\mu_1u \big),   \T (4.22)
$$
     where we set \f R_{0+}(\lambda) = R_0(\lambda) $. Since we have
$$
      |\mu_0^{-1}\mu_1u| \le c_1m_0^{-1}
                      (1 + |x|)^{-1-\epsilon+\delta-j\epsilon}
                       \big[(1 + |x|)^{-\delta+j\epsilon}|u|\big], \T (4.23)
$$
     it follows that 
$$
        \mu_0^{-1}\mu_1u \in  L_{2, 1-\delta+(j+1)\epsilon}(\RN).   \T (4.24)   
$$
     Noting that (4.19), \f 1/2 < \delta < 1/2 + \epsilon \A and 
     \f 0 < \epsilon < 1/2 $, we see that
$$
\left\{ \split
    \delta < 1/2 + \epsilon \Longrightarrow \ 1 - \delta + (j+1)\epsilon 
                                \ge  1 - \delta + \epsilon > 1/2, \\ *[5pt]
    - \delta + j\epsilon \le 0 \ {\rm and} \ 0 < \epsilon < 1/2 \\  *[5pt]
    \hspace{2.1cm} \Longrightarrow \  
                    1 - \delta + (j+1)\epsilon \le 1 + \epsilon < 3/2,
\endsplit \right.                                                 \T (4.25)
$$
     and hence we can apply Proposition 3.7 with \f \delta \A replaced by
     \f 1 - \delta + (j+1)\epsilon \A to obtain
$$
          \M \CD_r u \M_{-\delta+(j+1)\epsilon} \le
                       \lambda m_0^{-1}c_1C_j' \M u \M_{-\delta+j\epsilon}
                                                                   \T (4.26)
$$
     with \f C_j' = C_j'(\delta, \epsilon, m_0, M_0) \A for \f N \ge 3 $, and 
$$
          \M \CD_r u \M_{-\delta+(j+1)\epsilon,*} \le
                        \lambda m_0^{-1}c_1C_j'' \M u \M_{-\delta+j\epsilon}
                                                                   \T (4.27)
$$
      with \f C_j'' = C_j''(\delta, \epsilon, \lambda, m_0, M_0) \A for 
      \f N = 2 $.

        (II) Let \f R_0 > 1 $. 
     Set \f \beta = 2(-\delta + (j+1)\epsilon) + 1$, 
$$
     \xi(r) =
\left\{ \split
        r \ \ \ \  \hspace{1.5cm} \ \ \q (0 < r \le 1), \\ *[5pt]
        2^{-\beta}(1 + r)^{\beta} \ \ \ \q (1 < r \le R_0), \\ *[5pt]
        2^{-\beta}(1 + R_0)^{\beta} \ \q (r > R_0)
\endsplit \right.                                                 \T (4.28)
$$
     if \f N \ge 3 $, and
$$
     \xi(r) =
\left\{ \split
        {\dm 2^{-1} r^2 \ \  \hspace{1.5cm} \   \q (0 < r \le 1), } \\ *[5pt]
        2^{-\beta-1}(1 + r)^{\beta} \ \ \ \q (1 < r \le R_0), \\ *[5pt]
        2^{-\beta-1}(1 + R_0)^{\beta} \ \q (r > R_0)
\endsplit \right.                                                 \T (4.29)
$$
     for \f N = 2 $. It follows from (4.25) that \f 0 < \beta < 2 $, and 
     hence  \f \xi \A satisfies (a), (b), (c) or (c)' in Proposition 4.2. 
     Therefore we can apply (i) of Proposition 4.2 to obtain
$$
\split
    {\dm \int_{B_1} |u|^2 \,dx
          + \int_{B_{1R_0}} \beta 2^{-\beta}(1 + r)^{\beta-1}|u|^2 \,dx }
                                                                 \\ *[5pt]
     \hspace{1cm} {\dm  \le 2m_0^{-1}c_1 \int_{B_1} r(1 + r)^{-1-\epsilon}
                                           |u||\CD_ru|\,dx } \\ *[5pt]
     \hspace{2cm} {\dm + 2^{-\beta+1}m_0^{-1}c_1 \int_{B_{1R}} 
                      (1 + r)^{\beta-1-\epsilon}|u||\CD_ru|\,dx }  \\ *[5pt]
     \hspace{3cm} {\dm + C 2^{-\beta}(1 + R_0)^{\beta}
                        \int_{S_R} |\CD_ru|^2 \,dS }
\endsplit                                                          \T (4.30)
$$
     for \f R > R_0 \A and \f N \ge 3 $, or
$$
\split
    {\dm \int_{B_1} r|u|^2 \,dx
         + \int_{B_{1R_0}} \beta 2^{-\beta-1}(1 + r)^{\beta-1}|u|^2 \,dx }
                                                                   \\ *[5pt]
        \hspace{1cm} {\dm  \le m_0^{-1}c_1 \int_{B_1} 
                   r^2(1 + r)^{-1-\epsilon} |u||\CD_ru|\,dx } \\ *[5pt]   
        \hspace{2cm} {\dm + 2^{-\beta}m_0^{-1}c_1 \int_{B_{1R}} 
                      (1 + r)^{\beta-1-\epsilon}|u||\CD_ru|\,dx }  \\ *[5pt]
        \hspace{2cm} {\dm + \,2^{-\beta+1}m_0^{-1}\lambda^{-1} \int_{B_{1R}} 
                        (1 + r)^{2(-\delta+j\epsilon)}|u|^2\,dx } \\  *[5pt]
        \hspace{3cm} {\dm + C2^{-\beta-1}(1 + R_0)^{\beta}
                     \int_{S_R} |\CD_ru|^2 \,dS }
\endsplit                                                       \T (4.31)
$$
     for \f R > R_0 \A and\f N = 2,$ where we have used the facts that 
$$
\left\{ \split
       {\dm (1 + R_0)^{\beta} \le (1 + R)^{\beta}   \q  (R \ge R_0), } 
                                                              \\ *[5pt]  
       {\dm \frac{\varphi}r - 2^{-1}\frac{\pa\varphi}{\pa r} = 0 
                                                   \q (0 < r < 1), }
\endsplit \right.                                                 \T (4.32)   
$$
     and, for \f r > 1 $,    
$$
\split
      {\dm r^{-2} \bigg( \frac{\varphi}r - 2^{-1}\frac{\pa\varphi}{\pa r} 
                                                      \bigg) } \\  *[5pt]
       \hspace{1cm} {\dm \le 2^{-\beta-2}2^{3}\bigg(1 - \frac{\beta}2\bigg)
                              (1 + r)^{\beta-3} }   \\  *[7pt]
       \hspace{1cm} {\dm \le 2^{-\beta+1} (1 + r)^{2(-\delta+j\epsilon)}. }
\endsplit                                                         \T (4.33)
$$
     Since
$$
            \beta -1 - \epsilon = (-\delta + j\epsilon)
                            + (-\delta + (j+1)\epsilon),          \T (4.34)
$$
     we have
$$
    \int_{B_{1R}} (1 + r)^{\beta-1-\epsilon}|u||\CD_ru|\,dx
       \le \M u \M_{-\delta+j\epsilon}\M \CD_ru \M_{-\delta + (j+1)\epsilon}. 
                                                                   \T (4.35)
$$
     Therefore we have for \f N \ge 3 $
$$
\split
        {\dm \M u \M_{-\delta+(j+1)\epsilon, B_{R_0}}^2 }   \\  *[5pt]
          \hspace{1cm} {\dm \le C_0c_1\M u \M_{-\delta+j\epsilon}
                      \M \CD_ru \M_{-\delta+(j+1)\epsilon}
                     + \widetilde{C}(1 + R_0)^{\beta}
                        \int_{S_R} |\CD_ru|^2 \,dS, }    
\endsplit                                                           \T (4.36)
$$
     and for \f N = 2 $
$$
\split
        {\dm \M u \M_{-\delta+(j+1)\epsilon,B_{R_0}}^2 }   \\  *[5pt]
          \hspace{1cm} {\dm \le C_0' \bigg( c_1 \M u \M_{-\delta+j\epsilon}
                   \M \CD_ru \M_{-\delta+(j+1)\epsilon,*} 
                   + \lambda^{-1} \M u \M_{-\delta+j\epsilon}^2 \bigg) }   \\  *[5pt]
          \hspace{1cm} {\dm + \widetilde{C}'(1 + R_0)^{-\delta+j\epsilon}
                        \int_{S_R} |\CD_ru|^2 \,dS, }    
\endsplit                                                          \T (4.37)
$$
     where \f C_0 $, \f C_0' $, \f \widetilde{C} $, and \f \widetilde{C}' \A
     depend on \f j, \delta, \epsilon, m_0 $, and \f M_0 $. 
     Note that we have from (iii) of (4.2)
$$
           \liminf_{R\to\infty} \int_{S_R}  |\CD_ru|^2 \,dS = 0     \T (4.38)
$$
     Let \f R \to \infty \A along an appropriate sequence so that the last 
     terms of the right-hand sides of (4.36) and (4.37)
     tends to \f 0 $. Therefore, noting that \f R_0 > 1 \A is arbitrary, 
     and using (4.26) and (4.27), too, we obtain (4.21), which completes 
     the proof. \ \ \Box
\BP      
        
        {\bf Proposition 4.4.} {\it Assume that} \ Assumptions 2.1 {\it and}
     \ 2.2 {\it hold. Suppose that \f \mu_1 \A is long-range. 
     Suppose that \f u \in \widetilde{V}_p^{(+)}(H) \A [\, or \f u \in 
     \widetilde{V}_p^{(-)}(H) \,$] with an extended eigenvalue 
     \f \lambda \in \widetilde{\sigma}_p^{(+)}(H) \A [\, or 
     \f \lambda \in \widetilde{\sigma}_p^{(-)}(H) \, $]. 
     Let \f \epsilon' = \epsilon/2$. Suppose that
$$
\left\{ \split
              u \in L_{2, -\delta+j\epsilon'}(\RN),  \\ *[5pt]    
              - \delta + j\epsilon'  \le 0   
\endsplit \right.                                                   \T (4.39)
$$
     with a nonnegative integer \f j $. Then we have 
$$
       u \in L_{2, -\delta+(j+1)\epsilon'}(\RN),                    \T (4.40)
$$
     and
$$
       \M u \M_{-\delta+(j+1)\epsilon'} 
             \le \sqrt{c_1}C_j^{(N)}\M u \M_{-\delta+j\epsilon'}, \T (4.41) 
$$ 
     where \f C_j^{(N)} = C_j^{(N)}(\delta, \epsilon, m_0, M_0) \A for 
     \f N \ge 3 \A and  
     \f C_j^{(2)} = C_j^{(2)}(\lambda, \delta, \epsilon, m_0, M_0) $.}
\BP   

        Proof. \ \ We shall prove (4.40) and (4.41) in the case that 
     \f \lambda \in \widetilde{V}_p^{(+)}(H) $. We set 
     \f \CD_r^{(+)}u = \CD_ru $. Let \f R_0 > 1 $.
     Set \f \beta = 2(-\delta + (j+1)\epsilon') + 1$, and let \f \xi(r) \A 
     be given by (4.28) (with \f \epsilon \A replaced by \f \epsilon' $). 
     Here we should note that \f 0 < \beta < 3/2 $, and hence our 
     \f \xi(r) \A satisfies the conditions (a), (b) and (c) of 
     Proposition 4.2. Then we  have from (ii) of Proposition 4.2
$$
\split
     {\dm \int_{B_1} |u|^2 \, dx 
         + \int_{B_{1R_0}} \beta 2^{-\beta}(1+r)^{\beta-1} |u|^2 \, dx } 
                                                                  \\ *[6pt]
     \hspace{1cm} {\dm \le m_0^{-1}c_1 \int_{B_1} \big\{ (1+r)^{-\epsilon}
                      + r(1+r)^{-1-\epsilon} \big\} |u|^2 \, dx }  \\ *[6pt]
     \hspace{1.5cm} {\dm +  m_0^{-1}c_1(\beta+1) 2^{-\beta}
           \int_{B_{1R_0}} (1+r)^{\beta-1-\epsilon} |u|^2 \, dx } \\ *[6pt] 
     \hspace{1.5cm} {\dm +  m_0^{-1}c_1 2^{-\beta}
                    \int_{B_{1R_0}} (1+r)^{\beta}(1+r)^{-1-\epsilon} 
                                                     |u|^2 \, dx } \\ *[6pt]
     \hspace{1.5cm} {\dm + C2^{-\beta} (1+R_0)^{\beta}
                                    \int_{S_R} |\CD_ru|^2 \, dS, } \\ *[6pt] 
     \hspace{1cm} {\dm \le \widetilde{C} 
          \int_{B_{R}} (1+r)^{2(-\delta+j\epsilon')} |u|^2 \, dx } \\ *[6pt] 
     \hspace{1.5cm} {\dm + C2^{-\beta} (1+R_0)^{\beta}
                                    \int_{S_R} |\CD_ru|^2 \, dS, } 
\endsplit                                                          \T (4.42)
$$
     where \f R_0 < R $, 
     \f \widetilde{C} = \widetilde{C}(j, \delta, \epsilon, m_0) $, and
     we should note that \f \beta - 1- \epsilon = 2(-\delta+j\epsilon' $). The 
     inequality (4.41) follows from (4.42). The case that \f N = 2 \A can be 
     treated in quite a similar way, which completes the proof.  \ \ \Box 
\BP

        Now we are in a position to show that 
     \f V_p(H) = \widetilde{V}_p^{(\pm)}(H) $. Let
$$
\left\{ \split
         j_0 = \min \{ j \in \bN : -\delta + j\epsilon > 0 \, \},  \\ *[5pt]
         \delta_0 = -\delta + j_0\epsilon
\endsplit \right.                                                  \T (4.43)
$$
     if \f \mu_1 \A is short-range, and let
$$
\left\{ \split
       j_0 = \min \{ j \in \bN : -\delta + j\epsilon' > 0 \, \},  \\ *[5pt]
       \delta_0 = -\delta + j_0\epsilon'
\endsplit \right.                                                  \T (4.44)
$$
      if \f \mu_1 \A is long-range.
\BP
     
        {\bf Theorem 4.5.} {\it Let} \ Assumptions 2.1 {\it and} \ 2.2  
     {\it be satisfied.}  
\SP

        (i) \ {\it Then we have
$$
        \widetilde{V}_p^{(\pm)}(H) \subset H_{\delta_0}^2(\RN),    \T (4.45)
$$
     where \f \delta_0 \A is given by} \ (4.43) or (4.44), {\it and hence}
$$
\left\{ \split
        {\dm V_p(H) = \widetilde{V}_p^{(+)}(H) = \widetilde{V}_p^{(-)}(H), } 
                                                                    \\ *[5pt]
        {\dm \sigma_p(H) = \widetilde{\sigma}_p^{(+)}(H) 
                                    = \widetilde{\sigma}_p^{(-)}(H). }
\endsplit \right.                                                   \T (4.46)                                
$$
\SP

        (ii) \ {\it Let \f \mu_1 \A be short-range. Let \f u \in V_p(H) \A 
     associated with \f \lambda \in \sigma_p(H) $. Then, for each
     \f N \ge 2 $, there exists a positive constant \f C^{(N)} \A such that
$$
\left\{ \split
       {\dm \M u \M_{\delta_0} 
             \le C^{(N)}\big(c_1\sqrt{\lambda}\big)^{j_0}\M u \M_{-\delta} 
                                             \q (N \ge 3),} \\ *[5pt]
       {\dm \M u \M_{\delta_0} 
             \le C^{(2)}c_1^{j_0}\M u \M_{-\delta} 
                                             \q (N = 2),  }           
\endsplit \right.                                                 \T (4.47)
$$
     where \f j_0 \A is given by} \ (4.43), {\it and
     \f C^{(N)} = C^{(N)}(\delta, \epsilon, m_0, M_0) \A for \f N \ge 3 \A
     and  \f C^{(2)} = C^{(2)}(\lambda, \delta, \epsilon, m_0, M_0) $. 
     Further, for \f N \ge 3 $, we have}
$$
     \sigma_p(H) \subset [c_1^{-2}(C^{(N)})^{-2/j_0}, \ \infty).  \T (4.48)
$$
\SP

        (iii) \ {\it Let \f \mu_1 \A be long-range. Then, for each 
     \f N \ge 2 $, there exists a positive constant 
     \f C^{(N)} = C^{(N)}(\delta, \epsilon, m_0, M_0) \ (N \ge 3), \, = 
     C ^{(2)}(\lambda, \delta, \epsilon, m_0, M_0) \ (N = 2) \A such that
$$
        \M u \M_{\delta_0} \le C^{(N)}c_1^{j_0/2}\M u \M_{-\delta}, \T (4.49)
$$
     where \f u \in V_p(H) $, and \f j_0 \A and \f \delta_0 \A are given by} 
     (4.44). 
\BP

        Proof. \ \ Using Propositions 4.3 and 4.4 repeatedly, we obtain 
$$
         \widetilde{V}_p^{(\pm)}(H) \subset L_{2,\delta_0}(\RN),    \T (4.50)
$$
     and the inequalities (4.47) and (4.49), where 
     \f C^{(N)} = C_0^{(N)}C_1^{(N)} \cdots C_{j_0}^{(N)} \A and
     \f j_0 \A and \f \delta_0 \A are in (4.43) or (4.44). Let 
     \f u \in \widetilde{V}_p^{(\pm)}(H) \A associated with 
     \f \lambda \in \widetilde{\sigma}_p^{(\pm)}(H) $. Then it follows from 
     the equation \f -\,\Delta u - \lambda\mu u = 0 \A that 
     \f u, \ \Delta u \in L_{2,\delta_0}(\RN) \A which implies that 
     \f u \in H_{\delta_0}^2(\RN) $. Thus we have proved (4.45). 
     Let \f N \ge 3 \A and let \f \mu_1 \A is short-range. Since we have from 
     the first inequality of (4.47)
$$
         \M u \M_{-\delta} \le \M u \M_{\delta_0} 
             \le C^{(N)}\big(c_1\sqrt{\lambda}\big)^{j_0}\M u \M_{-\delta},
                                                                   \T (4.51)
$$
     or
$$
     (1 - C^{(N)}\big(c_1\sqrt{\lambda}\big)^{j_0})\M u \M_{-\delta} \le 0,
                                                                   \T (4.52)
$$
     whence (4.48) follows. This completes the proof. \ \ \Box
\BP             

        {\bf Theorem 4.6.}  {\it Let} \ Assumptions 2.1 {\it and} \ 2.2 
     {\it be satisfied. Let \f \sigma_p(H) \A be as above.}
\SP

        (i) \ {\it Then the multiplicity of each 
     \f \lambda \in \sigma_p(H) \A is finite.}
\SP

        (ii) {\it \f \sigma_p(H) \A does not have any accumulation point
     except \f \lambda = 0 \A and \f \lambda = \infty $. If \f N \ge 3 $, 
     then the only possible accumulation point of \f \sigma_p(H) \A is 
     \f \lambda = \infty $.}
\BP

        Proof. \ \ Suppose that \f \sigma_p(H) \A has an accumulation 
     point \f \lambda_0 \in (0, \ \infty) $. Then there exist infinite 
     sequences \f \{ \lambda_n \} \subset \sigma_p(H) \A and 
     \f  \{ u_n \} \subset V_p(H) \A such that
$$
\left\{ \split
          \lambda_n \to \lambda_0 \q (n \to \infty), \\ *[5pt]
         (u_m, \ u_n)_X = \delta_{mn} \q (m, n \in \bN), \\  *[5pt]
         -\,\mu(x)^{-1} \Delta u_n - \lambda_n u_n = 0 \q (n \in \bN),
\endsplit \right.                                                 \T (4.53)
$$
     where \f \delta_{mn} \A is Kronecker's delta. Since
$$
\split
       \M \gr u_n \M_0^2 = \lambda_n(u_n, u_n)_X = \lambda_n \M u_n \M_X^2 
                                                                   \\ *[5pt]
       \hspace{1.5cm}   = \lambda_n \le \sup_n \lambda_n < \infty,    
\endsplit                                                         \T (4.54)
$$
     we can apply the Rellich selection theorem to choose a
     subsequence \f \{ u_{n_m} \} \A which converges in 
     \f L_2(\RN)_{{\rm loc}} \A as \f m \to \infty $. Let 
     \f u_{0} \in L_2(\RN)_{{\rm loc}} \A be the limit function. On the 
     other hand, in view of Theorem 4.5, there exists a positive constant 
     \f C \A such that, for any \f s > 0 $,
$$
\split
     \M u_{n_m} \M_{0,E_s} 
            \le (1 + s)^{-\delta_0} \M u_{n_m} \M_{\delta_0,E_s}  
               \le (1 + s)^{-\delta_0} \M u_{n_m} \M_{\delta_0}  \\
     \hspace{1.8cm}       \le C(1 + s)^{-\delta_0} \M u_{n_m} \M_{-\delta} 
                                 \le C(1 + s)^{-\delta_0} \M u_{n_m} \M_{0}, 
\endsplit                                                    \T (4.55)
$$
     and hence 
$$
       \M u_{n_m} \M_{0,E_s} 
           \le \frac{C}{\sqrt{\widetilde{m_0}}}(1 + s)^{-\delta_0} \M u_{n_m} \M_{X}
           \le \frac{C}{\sqrt{\widetilde{m_0}}}(1 + s)^{-\delta_0},
                                                            \T (4.56)
$$
     where \f \delta_0 \A is given by (4.43) or (4.44). Therefore 
     \f u_{n_m} \A is small at infinity uniformly for \f m \in \bN $. 
     Thus it follows that \f u_{n_m} \A converges to \f u_0 \A in \f X \A 
     and \f \M u_0 \M_X = 1 $. Noting that  \f \{ u_{n_m} \} \A is an 
     orthonormal system in \f X $, we have
$$
     0 = \lim_{n\to\infty} (u_{n_m}, \ u_{n_{m+1}})_X =  \M u_0 \M_{X}^2 = 1,
                                                                   \T (4.57)
$$
     which is a contradiction. Therefore \f \sigma_p(H) \A is discrete in 
     \f (0, \ \infty) $. If \f N \ge3 $, (ii) of Theorem 4.5 implies that 
     \f \lambda = 0 \A cannot be an accumulation point of \f \sigma_p(H) $. 
     This completes the proof. \ \ \Box
\BP

        Consider the following additional condition on \f \mu_1(x) $:
\BP

        {\bf Assumption 4.7.} (i) The function \f \mu_1 \A is measurable 
     such that
$$
          \mu(x) \ge N\mu_1(x) + \lambda_0(|x||\mu_1(x)|)^2 \q ({\rm a.e. \ }
                                                     x \in \RN).  \T (4.58)
$$
     with \f \lambda_0 > 0 $.

        (ii) The function \f \mu_1 \A is differentiable
     and \f \mu_1 \A satisfies
$$
           \mu(x) + |x|\frac{\pa\mu_1}{\pa |x|} \ge 0 
                          \q ({\rm a. e. \ } x \in \RN).    \T (4.59)
$$
\BP
 
        The following theorem gives sufficient conditions that 
     the absence of \f \sigma_p(H) = 0 \A on some interval or whole positive
     half line (cf. Roach-Zhang\,[\RZ], the proof of Theorem 3.1). 
\BP
            
        {\bf Theorem 4.8.} {\it Suppose that 
     \f \mu(x) = \mu_0(x) + \mu_1(x) $, \f \mu_0 \A satisfies} 
     \ Assumptions 2.1, {\it and \f \mu \A satisfies} \ (2.12) {\it of} \ 
     Assumptions 2.2. {\it Suppose that} \ (i) {\it or} \ (ii)
     {\it of} \ Assumption 4.7 {\it hold. Then we have}
$$
\left\{ \split
       {\dm \sigma_p(H) \cap [0, \lambda_0] = \emptyset
                                  \qq (\mbox{if \ (i) \ holds}),} \\ *[5pt]
       {\dm \sigma_p(H) = \emptyset  \qq \ \ \ \ \ \ \ \ \ \ \
                                     (\mbox{if \ (ii) \ holds}).}                                      
\endsplit \right.                                               \T (4.60)
$$
\BP

        Proof. \ \ (I) Let \f u \in H^2(\RN) \A satisfy the homogeneous
     equation $ -\,\Delta u - \lambda\mu(x)u = 0 $ with \f \lambda > 0 $. 
     We have only to show that \f u \equiv 0 $.
     We are going to multiply both sides of the equation by 
     \f 2r(\pa_r\overline{u}) + (N - 1)\overline{u} $, integrate over 
     \f B_R $, \f R > 0 $, and take the real part. 

        (II) Using the identity
$$
         2{\rm Re}\,[(\Delta u)r(\pa_r\overline{u})]
                 = {\rm div}\,[2{\rm Re}\{r(\pa_r\overline{u})\gr u\}
                                            - |\gr u|^2x]
                                       + (N - 2)|\gr u|^2     \T (4.61)
$$
     (Roach-Zhang\, [\RZ], (3.4) with \f h(r) \equiv 1 $) and the divergence
     theorem, we have
$$
\split
        {\dm 2{\rm Re}\,\int_{B_R} (-\Delta u) r(\pa_r\overline{u}) \, dx }
                                                                   \\ *[6pt]
        \hspace{1.5cm} {\dm  = - \int_{B_R} (N - 2)|\gr u|^2 \, dx
                          - R \int_{S_R} (2|\pa_r\overline{u}|^2 
                                            - |\gr u|^2) \, dS,}
\endsplit                                                     \T (4.62)
$$
     where \f \pa_r v = \pa v/\pa r $, \f r = |x| $. Since it is easy to see 
     that
$$
\split
  {\dm{\rm Re}\,\int_{B_R} (-\Delta u)(N - 1) \overline{u} \, dx }  \\ *[6pt]
   \hspace{1.5cm} {\dm                                              
             = \int_{B_R} (N - 1)|\gr u|^2 \, dx -
                (N - 1) \int_{S_R} {\rm Re}\,[(\pa_ru)\overline{u}] \, dS,}
\endsplit                                                         \T (4.63)
$$
     it follows that 
$$
\split
        {\dm {\rm Re}\,\int_{B_R} (-\Delta u) \{ 2r(\pa_r\overline{u})
                                       + (N - 1) \overline{u} \} \, dx }
                                                                   \\ *[6pt]    
        \hspace{1.5cm} {\dm   = \int_{B_R} |\gr u|^2 \, dx
                       -  R \int_{S_R} (2|\pa_r\overline{u}|^2 - |\gr u|^2
                      + \frac{N-1}R{\rm Re}\,[(\pa_ru)\overline{u}]) \, dS.}
\endsplit                                                         \T (4.64)
$$

        (III) Suppose that (ii) of Assumption 4.7 holds. By the use of the 
     integration by parts, we have
$$
\split
    {\dm 2{\rm Re}\,\int_{B_R} (-\lambda\mu u)r(\pa_r\overline{u}) \, dx } 
                                                                  \\ *[6pt]    
    \hspace{1cm} {\dm  = - \, \int_{B_R} (\lambda\mu)r(\pa_r|u|^2) \, dx } 
                                                                 \\ *[6pt]    
    \hspace{1cm} {\dm = \lambda \int_{B_R} (N\mu + r(\pa_r\mu_1))|u|^2 \, dx} 
                                                                \\ *[6pt]
    \hspace{2cm} {\dm - \lambda \sum_{\ell\in L} 
                          \int_{\pa\Omega_{\ell}\cap B_R}
                            \mu_0(x\cdot n^{(\ell)}) |u|^2 \, dS 
                         - \, \lambda R\,\int_{S_R} \mu |u|^2 \, dS, } 
\endsplit                                                        \T (4.65)
$$
    where we should note that \f \mu_0 \A does not appear in the first term
    of the right-hand side since it is constant on each \f \Omega_{\ell} $, 
    and \f \mu_1 \A does not appear in the second term of the right-hand 
    side since it is continuous on \f \RN $. Also we should note that
$$
\split
    {\dm - \lambda \sum_{\ell\in L} \int_{\pa\Omega_{\ell}\cap B_R}
                       \mu_0(x\cdot n^{(\ell)}) |u|^2 \, dS }    \\ *[6pt]  
    \hspace{1cm} {\dm  = \lambda \sum_{\ell\in L} 
                    \int_{S_{\ell}^{(+)}\cap B_R}
         (\nu_{\ell+1} - \nu_{\ell}) (x\cdot n^{(\ell)}) |u|^2 \, dS \ge 0 }
\endsplit                                                        \T (4.66)
$$
     since the integrand is nonnegative by (2.11). Thus, 
$$
\split
    {\dm {\rm Re}\,\int_{B_R} (-\lambda\mu u)[ 2r(\pa_r\overline{u}) 
                                 + (N -1)\overline{u}]\, dx  }  \\ *[6pt]        
    \hspace{1cm} {\dm = \lambda \int_{B_R} (\mu + r(\pa_r\mu_1))|u|^2 \, dx } 
                                                                 \\ *[6pt]
    \hspace{1.5cm} {\dm + \lambda \sum_{\ell\in L} 
                   \int_{S_{\ell}^{(+)}\cap B_R}
                (\nu_{\ell+1} - \nu_{\ell}) (x\cdot n^{(\ell)}) |u|^2 \, dS
                         - \, \lambda R\,\int_{S_R} \mu |u|^2 \, dS } 
\endsplit                                                         \T (4.67)
$$

        (IV) It follows from (4.64) and (4.67) that
$$
\split
    {\dm 0 = {\rm Re}\,\int_{B_R} (-\Delta u - \lambda\mu u) 
              \{ 2r(\pa_r\overline{u}) + (N - 1) \overline{u} \} \, dx }
                                                                   \\ *[6pt]                                           
    \hspace{.3cm} {\dm = \int_{B_R} \big(|\gr u|^2 
                          + \lambda(\mu + r(\pa_r\mu_1))|u|^2 \big)\, dx }    
                                                              \\ *[6pt]  
    \hspace{.8cm} {\dm + \lambda \sum_{\ell\in L} 
                   \int_{S_{\ell}^{(+)}\cap B_R}
             (\nu_{\ell+1} - \nu_{\ell}) (x\cdot n^{(\ell)}) |u|^2 \, dS }  
                                                               \\ *[6pt] 
    \hspace{.8cm} {\dm + R \int_{S_R} (|\gr u|^2 - 2|\pa_r\overline{u}|^2 
                      - \frac{N-1}R{\rm Re}\,[(\pa_ru)\overline{u}]
                            - \lambda\mu |u|^2) \, dS.   }
\endsplit                                                       \T (4.68)
$$
     Since \f 2|\pa_r\overline{u}|^2 + |\gr u|^2 + \lambda\mu |u|^2 
     + (N-1)|\pa_ru||u| \A is integrable on \f \RN $, we see that the third 
     term of the right-hand side goes to \f 0 \A as \f R \to \infty \A 
     along an appropriate sequence. Therefore it follows from (4.68) that
$$
\split
    {\dm 0 = \int_{\RN} \big(|\gr u|^2 
               + \lambda(\mu + r(\pa_r\mu_1))|u|^2 \big)\, dx }    \\ *[6pt]  
    \hspace{.8cm} {\dm + \lambda \sum_{\ell\in L} \int_{S_{\ell}^{(+)}}
             (\nu_{\ell+1} - \nu_{\ell}) (x\cdot n^{(\ell)}) |u|^2 \, dS. } 
\endsplit                                                        \T (4.69)
$$
     Noting that all the integrands in the right-hand side are nonnegative, 
     we have \f \gr u = 0 \A a. e., and hence \f u \equiv 0 \A since 
     \f u \in H^2(\RN) $.
    
        (V) Suppose that (i) of Assumption 4.7 holds. By using partial 
     integration only for the term containing \f \mu_0 $, we obtain
$$
\split
    {\dm {\rm Re}\,\int_{B_R} (-\lambda\mu u) \big[ 2r(\pa_r\overline{u}) 
                                         + (N - 1)\overline{u} \big]\, dx } 
                                                                  \\ *[6pt]    
    \hspace{1cm} {\dm  = - \lambda\, 
                  \int_{B_R} \big[ \mu_0 r(\pa_r|u|^2) 
                                 + (N - 1)\mu_0|u|^2 \big]\, dx } \\ *[6pt]
     \hspace{2cm} {\dm - \lambda\, 
                      \int_{B_R} \mu_1\big[ 2r{\rm Re}(u(\pa_r\overline{u})) 
                                         + (N - 1)|u|^2 \big]\, dx } 
                                                                 \\ *[6pt]    
    \hspace{1cm} {\dm = \lambda \int_{B_R} \mu_0|u|^2 \, dx  } 
                                                                  \\ *[6pt]
    \hspace{2cm} {\dm - \lambda \sum_{\ell\in L} 
                        \int_{\pa\Omega_{\ell}\cap B_R}
                            \mu_0(x\cdot n^{(\ell)}) |u|^2 \, dS 
                  - \, \lambda R\,\int_{S_R} \mu_0 |u|^2 \, dS, } \\ *[6pt]
    \hspace{2cm} {\dm  - \lambda\, 
                    \int_{B_R} \mu_1\big[ 2r{\rm Re}(u(\pa_r\overline{u})) 
                                         + (N - 1)|u|^2 \big]\, dx }
\endsplit                                                        \T (4.70)
$$
     Let \f h(x) \A be a positive function to be specified later. Since 
     we have
$$
     2|r\mu_1u(\pa_r\overline{u})| \le |r\mu_1| \big( h|u|^2 + 
                                        \frac{|\gr u|^2}h \big), \T (4.71)
$$
    it follows from (4.70) that
$$
\split
    {\dm {\rm Re}\,\int_{B_R} (-\lambda\mu u) \big[ 2r(\pa_r\overline{u}) 
                                       + (N - 1)\overline{u} \big]\, dx } 
                                                               \\ *[6pt]    
    \hspace{1cm} {\dm  \ge - \lambda\, 
                  \int_{B_R} \frac{|r\mu_1|}{h}|\gr u|^2 \, dx } \\ *[6pt]
    \hspace{2cm} {\dm + \lambda \, \int_{B_R} \big[ \mu_0 - (N - 1)\mu_1 
                       - h|r\mu_1| \big]|u|^2 \, dx } \\ *[6pt]
    \hspace{2cm} {\dm - \lambda \sum_{\ell\in L} 
                         \int_{\pa\Omega_{\ell}\cap B_R}
                        \mu_0(x\cdot n^{(\ell)}) |u|^2 \, dS 
                - \, \lambda R\,\int_{S_R} \mu_0 |u|^2 \, dS. } \\ *[6pt]
\endsplit                                                       \T (4.72)
$$
     Thus (4.64) and (4.72) are combined to give
     
\newpage
     
$$
\split
    {\dm 0 = {\rm Re}\,\int_{B_R} (-\Delta u -\lambda\mu u) 
              \big[ 2r(\pa_r\overline{u}) + (N - 1)\overline{u} \big]\, dx } 
                                                                 \\ *[6pt]    
    \hspace{.4cm} {\dm  \ge \int_{B_R} 
              \big( 1 - \frac{\lambda|r\mu_1|}{h} \big)|\gr u|^2 \, dx } 
                                                                 \\ *[6pt]
    \hspace{1.0cm} {\dm + \lambda \, \int_{B_R} \big[ \mu - N\mu_1 
                       - h|r\mu_1| \big]|u|^2 \, dx } \\ *[6pt]
    \hspace{1.5cm} {\dm - \lambda \sum_{\ell\in L} 
                                 \int_{\pa\Omega_{\ell}\cap B_R}
                            \mu_0(x\cdot n^{(\ell)}) |u|^2 \, dS 
                 - \, \lambda R\,\int_{S_R} \mu_0 |u|^2 \, dS } \\ *[6pt]
    \hspace{1.5cm} {\dm -  R \int_{S_R} (2|\pa_r\overline{u}|^2 - |\gr u|^2
                      + \frac{N-1}R{\rm Re}\,[(\pa_ru)\overline{u}]) \, dS.}              
\endsplit                                                         \T (4.73)
$$
     Then, letting \f R \to \infty \A along an appropriate sequence in 
     (4.73), we obtain
$$
\split 
      {\dm 0 \ge \int_{\RN} \big( 1 - \frac{\lambda|r\mu_1|}{h} \big)
                                               |\gr u|^2 \, dx} \\ *[6pt]

      \hspace{1.1cm} {\dm + \lambda \, \int_{\RN} \big[ \mu - N\mu_1 
                       - h|r\mu_1| \big]|u|^2 \, dx },
\endsplit                                                        \T (4.74)
$$
     where we have used (4.66), too. Let \f \lambda \in [0, \lambda_0) \A
     be an eigenvalue of \f H \A with its eigenfunction \f u $.
     Set \f \eta = \lambda/\lambda_0 \in [0, 1) \A and
$$
       h(x) =
\left\{ \split 
       {\dm 1 \ \ \ \ \ \ \ \ \ \ \ \ 
                           \qq ({\rm if } \ \mu_1(x) = 0 ),} \\ *[6pt]
       {\dm \lambda_0|r\mu_1(x)| \qq ({\rm if } \ \mu_1(x) \ne 0 ).}
\endsplit \right.                                               \T (4.75)
$$
     Then we have 
$$
      {\dm 1 - \frac{\lambda|r\mu_1(x)|}{h(x)} = } 
\left\{ \split
         {\dm 1 \ \ \ \ \    \qq ({\rm if } \ \mu_1(x) = 0 ),} \\ *[6pt]
         {\dm 1 - \eta  \qq ({\rm if } \ \mu_1(x) \ne 0 ),}
\endsplit \right.                                              \T (4.76)
$$
     and hence, by using (4.58) 
$$
    {\dm \mu(x) - N\mu_1(x) - h|r\mu_1(x)| = } 
\left\{ \split
    {\dm \mu_0 > 0 \ \ \ \ \ \ \ \ \        
                     \qq ({\rm if } \ \mu_1(x) = 0 ),} \\ *[6pt]
    {\dm \mu(x) - N\mu_1(x) - \lambda_0(r\mu_1(x))^2 \ge 0} \\ *[6pt]
    {\dm       \ \ \ \ \ \ \ \ \ \ \ \ \ \ \ \ \ \ \ \ \ \ \  
                 \q   ({\rm if } \ \mu_1(x) \ne 0 ).}
\endsplit \right.                                                 \T (4.77)
$$       
     Therefore, we have from (4.74) 
$$
         0 \ge (1 - \eta) \int_{\RN} |\gr u|^2 \, dx,              \T (4.78)
$$
     i.e., \f \gr u \equiv 0 \A or \f u \A is identically zero almost 
     everywhere. This completes the proof.  \ \ \Box

$$
\ \ \ \
$$

\newpage

        {\bf 5. \ The limiting absorption principle for \f H $}              
\BP

        Throughout this section we assume that \f \delta \A satisfies
$$
             \frac12 < \delta \le \frac{1}2 + \frac{\epsilon}4,      \T (5.1)
$$
     where \f \epsilon \A is as in (2.13) or (2.14). Let \f u \in X \A be 
     given by
$$
\left\{ \split
          u = R(z)f,  \\ *[4pt]
          z = \lambda + i\eta  \qq (\lambda \ge 0, \eta \ne 0),  \\ *[4pt]
          f \in L_{2, \delta}(\RN),
\endsplit \right.                                            \T (5.2)
$$ 
     where \f R(z) = (H - z)^{-1} $. Then \f u \A satisfies the inhomogeneous
     equation  \f (-\mu^{-1}\Delta - z)u = f \A which is equivalent to
$$
                   (-\mu_0^{-1}\Delta - z)u = g \q
                            (g = \mu_0^{-1}(\mu f + z\mu_1u))   \T (5.3)
$$
     with \f k = \sqrt{z\mu_0} $. Let \f \mu_1 \A be short-range. Then, since
$$
            u \in L_{2, -\delta}(\RN) \Longrightarrow \ 
\left\{ \split
                    \mu_1u \in L_{2, \delta}(\RN),  \\ *[4pt]   
                    \M \mu_1u \M_{\delta} \le c_1 \M u \M_{-\delta},
\endsplit \right.                                                    \T (5.4)
$$
     we see that \f g \in L_{2, -\delta}(\RN) $. In the case that 
     \f \mu_1 \A is long-range, the inequality
$$
          |\eta|\M u \M_{1,0}^2 \le C(|f|, \, |u|)_0               \T (5.5)
$$
     will be useful, where \f C = C(\mu) $, \f \M \ \M_{1,0} \A is the norm 
     of \f H^1(\RN) $, and \f ( \ , \ )_0 \A is the inner product of 
     \f L_2(\RN) $. For the proof of (5.5), see, e.g., Eidus\,[\Ei], [\Sd],
     Lemma 2.1. Then, by a direct application of Theorem 3.5 
     to our case, we can evaluate the radiation condition term \f \CD u $.
\BP

        {\bf Theorem 5.1.} {\it Suppose that} \ Assumptions 2.1 {\it and} \
     2.2 {\it hold. Let \f \delta \A be as in} \ (5.1). {\it Let 
     \f 0 < c < d < \infty \A and let \f J_{\pm}(c, d) \A be as in} \ 
     (3.18). {\it Let \f u \A be given by} \ (5.2) {\it with 
     \f z \in J_{+}(c, d) \cup J_{-}(c, d) $. Then there
     exists a positive constant \f C = C(\delta, c, d, m_0, M_0) \A such 
     that 
$$
        \M \CD u \M_{\delta-1} 
                  \le C \big( \M f \M_{\delta} + \M u \M_{-\delta} \big)
                                                                 \T (5.6)
$$
     for \f N \ge 3 $, and
$$
        \M \CD u \M_{\delta-1, *} 
                  \le C \big( \M f \M_{\delta} + \M u \M_{-\delta} \big)
                                                                 \T (5.7)
$$     
     for \f N = 2 $, where \f \M  \ \M_{t, *} \A is as in} \ (3.21).
\BP

        Proof. \ \ We can proceed as in the proof of Theorem
     3.5. Since \f f \A in the proof of Theorem 3.5 should be replaced by
     \f g = \mu_0^{-1}(\mu f + z\mu_1u)\A (see (5.3)), our additional task is to 
     prove, for any \f e > 0 $,       
$$
     I = {\rm Re} \int_{B_{rR}} z\varphi\mu_1(x)u\overline{\CD_ru}\, dx 
              \le C( \M f \M_{\delta}^2 + \M u \M_{-\delta}^2
                              + e\M \CD u \M_{\delta-1}^2 )        \T (5.8)
$$
     with \f C = C(e, \delta, c, d, m_0, M_0) $, where \f \varphi \A is as in 
     the proof of Theorem 3.5, and \f \M \CD u \M_{\delta-1} \A
     in (5.8) should read \f \M \CD u \M_{\delta-1,*} \A if \f N = 2 $. 
     Suppose that \f \mu_1 \A is short-range. Then (5.8) follows directly 
     from (5.4). Suppose that \f \mu_1 \A is long-range. Then we have from 
     the definition of \f \CD_ru \A ((6) of Notation 3.1) and partial 
     integration
$$
\split
   {\dm I = \lambda {\rm Re}\int_{B_{rR}} \varphi\mu_1(x)u\overline{\CD_ru}\, dx 
                 - \eta {\rm Im} \int_{B_{rR}} 
                           \varphi\mu_1(x)u\overline{\CD_ru}\, dx} \\ *[6pt]
   \hspace{.3cm} {\dm = - \frac{\lambda}2 
                                \int_{B_{rR}} \pa_r(\varphi\mu_1)|u|^2 \, dx 
                        + \lambda \int_{B_{rR}} b\varphi\mu_1|u|^2 \, dx } 
                                                                   \\ *[6pt]
   \hspace{2.3cm} {\dm   - \eta {\rm Im} \int_{B_{rR}} 
                           \varphi\mu_1(x)u\overline{\CD_ru}\, dx 
                                 + I_4(r) + I_5(R) } \\ *[6pt]
   \hspace{.3cm} {\dm = I_1 + I_2 + I_3 + I_4(r) + I_5(R), }
\endsplit                                                          \T (5.9)
$$                    
     where the terms \f I_4(r) \A and \f I_5(R) \A tend to zero as 
     \f r \to 0 \A and \f R \to \infty \A along appropriate sequences, 
     respectively. It follows from (5.1) and the definition of 
     \f \varphi \A ((3.22) or (3.23)) that \f \varphi\mu_1 \A is bounded on 
     \f \RN $, and hence \f I_2 \A and \f I_3 \A can be evaluated by using 
     (5.5). On the other hand, since
$$
          \pa_r(\varphi\mu_1) = O((1 + r)^{2\delta-2-\epsilon})
                                 = O((1 + r)^{-2\delta}),          \T (5.10)
$$
     the term \f I_1 \A is evaluated by \f \M u \M_{-\delta}^2 $, which 
     completes the proof.  \ \ \Box
\BP
                                  
        As in \S4, let \f \sigma_p(H) \A be the set of all eigenvalues
     of \f H \A which is a discrete set in \f (0, \, \infty) $ (Theorem 4.6).
     Let \f \lambda > 0 \A such that \f \lambda \not\in \sigma_p(H) $. Let 
     \f u \in H^2(\RN)_{{\rm loc}} \cap L_{2, -\delta}(\RN) \A be a 
     solution of the homogeneous equation 
     \f -\,\mu(x)^{-1}\Delta u - \lambda u = 0 \A with the radiation 
     condition \f \M \CD^{(+)}u \M_{\delta-1, E_1} < \infty \A or
     \f \M \CD^{(-)}u \M_{\delta-1, E_1} < \infty $. Then it follows from 
     Theorem 4.5 that \f u \in L_{2, \delta_0}(\RN) \A where \f \delta_0 \A
     is given by (4.43) or (4.44). Since \f \lambda \A is supposed not to be
     an eigenvalue, we have \f u \equiv 0 $. Therefore we can prove the 
     limiting absorption principle for 
     \f \lambda \in (0, \infty) \bs \sigma_p(H) \A by starting with Theorem
     5.1, proceeding as in \S5 $\sim$ \S7 of [\JS]. Let 
     \f D_{\pm} \subset \bC \A be given by (3.29). For \f \lambda > 0$, let 
$$
            R_{\pm}(\lambda) = \lim_{\eta\downarrow 0} R(\lambda \pm i\eta),   
                                                                    \T (5.11)
$$
     and extend the resolvent \f R(z) \A on \f D_{\pm} \A by
$$
     R(\lambda+i\eta) = 
\left\{ \split
    R(\lambda+i\eta) \qq (\lambda > 0, \, \eta > 0), \\  *[4pt]      
    R_{+}(\lambda) \ \ \qq \ \ (\lambda > 0, \, \eta = 0)    
\endsplit \right.                                               \T (5.12)
$$
     for \f z \in D_+ \A and                   
$$
     R(\lambda+i\eta) = 
\left\{ \split
    R(\lambda+i\eta) \qq (\lambda > 0, \, \eta < 0), \\ *[4pt]       
    R_{-}(\lambda) \ \ \qq \ \ (\lambda > 0, \, \eta = 0)    
\endsplit \right.                                               \T (5.13)
$$
     for \f z \in D_- $. Then we have

\newpage

        {\bf Theorem 5.2.} {\it Suppose that} \ Assumptions 2.1 {\it and} \
     2.2 {\it holds. Let \f \delta \A satisfy} \ (5.1). 
\SP

        (i) \ {\it Then the limits} \ (5.11) {\it is well-defined in
     \f \bB(L_{2,\delta}(\RN), \, H_{-\delta}^2(\RN)) \A for
     \f \lambda \in (0, \infty) \bs \sigma_p(H) $, and the extended
     resolvent \f R(z) \A is a 
     \f \bB(L_{2,\delta}(\RN), \, H_{-\delta}^2(\RN))$-valued 
     continuous function on each of \f D_+ \bs \sigma_p(H) \A and 
     \f D_- \bs \sigma_p(H) $. }
\SP

        (ii) \ {\it For any \f z \in D_+ \bs \sigma_p(H) \A 
     {\rm [}\,or \f D_- \bs \sigma_p(H) $\,{\rm ]}, 
     \f R(z) \A is a compact operator from \f L_{2,\delta}(\RN) \A 
     into \f H_{-\delta}^1(\RN) $.  }   
\SP

        (iii) \ {\it The selfadjoint operator \f H \A is absolutely 
     continuous on the interval \f [c, d] \A such that 
     \f 0 < c < d < \infty \A and 
$$
            [c, \, d] \cap \sigma_p(H) = \emptyset.            \T (5.14)
$$         
     The operator \f H \A has no singular continuous spectrum.}
\SP

        (iv) \ {\it For \f 0 < c < d < \infty \A satisfying} \ (5.14) 
     {\it there exists \f C = C(c, d, \delta, m_0, M_0) > 0 \A such 
     that, for 
     \f z \in \overline{J}_+(c, \, d) \cup \overline{J}_-(c, \, d) $,
$$   
\left\{ \split
      {\dm \int_{E_s}(1 + r)^{-2\delta} 
              \big( |\gr R(z)f|^2 + |k|^2|R(z)f|^2 \big) \, dx } \\  *[5pt]
      \hspace{3.5cm} 
          {\dm  \le C^2 (1 + s)^{-(2\delta-1)} \M f \M_{\delta}^2 }\\ *[4pt]
       \ \ \ \ \ \ \ \ \ \ \ \ \ \ \ \ \ \ \ \ \ \ \ \ \ \ \ \ \ \ \ \ \
         \ \ \ \ \ \ \ \ {\dm (s \ge 1, \ f \in L_{2,\delta}(\RN)), } 
                                                                 \\  *[4pt] 
     {\dm \M \CD R(z)f \M_{\delta-1} \le C \M f \M_{\delta} \q
             ( f \in L_{2,\delta}(\RN)), }     
\endsplit \right.                                              \T (5.15)
$$
     where, for \f \lambda \in D_+ \cap (0, \, \infty) \A 
     {\rm [}\,or \f D_- \cap (0, \, \infty)${\rm ]}, \f \CD u \A should 
     be interpreted as \f \CD^{(+)} \A {\rm [}\,or \f \CD^{(-)}\,${\rm ]},
     and \f \overline{J}_{\pm}(c, \, d) \A are given by} \ (3.30). 
\BP

\BP

        {\sc Willi J\"ager \hspace{.65in} University of Heidelberg, Germany}
\SP
      
        {\sc Yoshimi Saito \hspace{.5in} University of Alabama at Birmingham,
                                                                  U. S. A.} 
$$
\ \ \ 
$$ 
   
            
\cen{  {\bf References}  }           
\BP

\F   [1] S. Agmon, {\it Spectral properties of Schr\"{o}dinger operators and 
     scattering theory}, Ann. 
     
\F   \ \ \ \ Scuola Sup. Pisa {\bf 2} (1975), 151-218. 
\SP

\F   [2] M. Ben-Artzi, Y. Dermanjian and J.-C. Guillot, {\it Acoustic waves 
     in perturbed stratified 
     
\F   \ \ \ fluids: a spectral theory}, Commun. Partial Differential 
     Equations {\bf 14} (1989), 479-517.  
\SP

\F   [3] A. Boutet de Monvel-Berthier and D, Manda, {\it Spectral and 
     scattering theory for wave 
     
\F   \ \ \ propagation in perturbed stratified media}, Universit\"{a}t 
     Bielefeld, BiBos. Preprint, Nr. 
     
\F   \ \ \ \ 606/11/93.
\SP

\F   [4] S. DeBi\'{e}vre and D. W. Pravica, {\it Spectral analysis for 
     optical fibres and stratified fluids 
     
\F   \ \ \ \ I: The limiting absorption principle}, J. Functional Analysis 
     {\bf 98} (1991) 406-436.
\SP

\F   [5] S. DeBi\'{e}vre and D. W. Pravica, {\it Spectral analysis for 
     optical fibres and stratified fluids 
     
\F   \ \ \ \ II: Absence of eigenvalues}, Commun. Partial Differential 
     Equations {\bf 17}, (1992), 69-97.   
\SP

\F   [6] D. Eidus, {\it The limiting absorption and amplitude problems for 
     the diffraction problem 
     
\F   \ \ \ with two unbounded media}, Comm. Math. Phys. {\bf 107} 
     (1986), 29-38.
\SP

\F   [7] T. Ikebe and Y. Sait\={o}, {\it Limiting absorption method and 
     absolute continuity for the 
     
\F   \ \ \ \ Schr\"odinger operator}, J. Math. Kyoto Univ. {\bf 12} (1972), 
     513-542.
\SP
             
\F   [8] W. J\"{a}ger and Y. Sait\={o}, {\it The limiting absorption 
     principle for the reduced wave operator 
     
\F   \ \ \ \ with cylindrical discontinuity}. Seminar Notes. IWR (SFB359)
     94-74, University of Heidelberg. 1994.
\SP

\F   [9] W. J\"{a}ger and Y. Sait\={o}, {\it On the Spectrum of the
     reduced wave operator with cylindrical 
     
\F   \ \ \ \ discontinuity}, Forum Mathematicum {\bf 9} (1997), 29-60.
\SP
 
\F   [10] G. Roach and B. Zhang, {\it On Sommerfeld radiation conditions for 
     the diffraction prob-
     
\F   \ \ \ \ \, lem with two unbounded media}, Proc. Royal Soc. Edinburgh 
     {\bf 121A} (1992), 149-161.
\SP
 
\F   [11] Y. Sait\={o}, {\it The principle of limiting absorption for 
     second-order differential equations 
     
\F   \ \ \ \ \, with operator-valued coefficients}, Publ. RIMS, Kyoto Univ. 
     {\bf 7} (1972), 581-619.
\SP

\F   [12] Y. Sait\={o}, {\it The principle of limiting absorption for the 
     non-selfadjoint Schr\"odinger 
     
\F   \ \ \ \ \, operators in \f \RN \A (\f N \neq 2 $)}, Publ. RIMS, Kyoto Univ.
     {\bf 9} (1974), 397-428.
\SP

\F   [13] Y. Sait\={o}, {\it A remark on the limiting absorption principle 
     for the reduced wave equation 
     
\F   \ \ \ \ \, with two unbounded media,} Pacific J. Math. {\bf 136} (1989), 
     183-208.

\SP

\F   [14] R. Weder, {\it Absence of eigenvalues of the acoustic propagators 
     in deformed waveguides}, 
     
\F   \ \ \ \ \, Rocky Mountain J. Math. {\bf 18} (1988), 495-503.
\SP

\F   [15] R. Weder, Spectral and Scattering Theory for Wave Propagation 
     in Perturbed Strat-
     
\F   \ \ \ \ \, ified Media, Springer-Verlaga, Berlin, 1991.
\SP

\F   [16] C. Wilcox, Sound Propagation in Stratified Fluids, 
     Springer-Verlag, New York, 1984.
\SP

\F   [17] B. Zhang, {\it On radiation conditions for acoustic propagators in 
     perturbed stratified 
     
\F   \ \ \ \ \, fluids}. Preprint. 1994. 
          .

\end{document}